\documentclass[12pt]{article}
\usepackage[utf8]{inputenc}
\usepackage[margin=1in]{geometry}
\usepackage{graphicx}
\usepackage{amsmath,amssymb,amsthm,authblk,mathrsfs,subfigure}
\usepackage[acronym]{glossaries}
\usepackage{color}
\usepackage{tikz} 
\usepackage[english]{babel}
\usepackage{comment}
\usepackage[nottoc]{tocbibind} 
\usepackage[outdir=./]{epstopdf}

\usepackage{url}
\usepackage{thmtools}

\usepackage{setspace}

\newcommand{\be}{\begin{equation}}
\newcommand{\ee}{\end{equation}}

\newcommand {\R}{\mathbb R}
\newcommand {\bu}{\bar u}
\newcommand {\bx}{\bar x}
\newcommand {\cT}{\mathcal T}

\newcommand {\sgn}{\mbox{sgn}\,}

\theoremstyle{plain}

\makeatletter
\newcommand{\sq}[1]{\mathbin{\mathpalette\make@circled{#1}}} 
\newcommand{\make@circled}[2]{%
	\ooalign{$\m@th#1\smallbigcirc{#1}$\cr\hidewidth$\m@th#1#2$\hidewidth\cr}%
}
\newcommand{\smallbigcirc}[1]{%
	\vcenter{\hbox{\scalebox{1.2}{$\m@th#1\square$}}}%
}
\makeatother

\declaretheorem[name={Example},qed={\lower-0.3ex\hbox{$\square$}} ] {Example}

\newcommand{\st}{\, | \,}

\title{Maximizing average throughput in oscillatory biological synthesis systems: an optimal control approach}
\author[1]{\normalsize M. Ali Al-Radhawi}
\author[2]{\normalsize Michael Margaliot}
\author[1,3,*]{\normalsize  Eduardo D. Sontag}
\affil[1]{\normalsize Departments of Bioengineering, and Electrical and Computer Engineering, Northeastern University, Boston MA 02115.}
\affil[2]{\normalsize Department of Electrical Engineering-Systems, Tel Aviv University}
\affil[3]{\normalsize Laboratory of Systems Pharmacology, Program in Therapeutic Science, Harvard Medical School, Boston MA 02115}
\affil[*]{\normalsize Corresponding Author: \texttt{sontag@sontaglab.org} }
\date{\today}
\newtheorem{theorem}{Theorem} 
\newtheorem{lemma}[theorem]{Lemma}
\newtheorem{problem}{Problem}

\newtheorem{proposition}[theorem]{Proposition} 
\newtheorem{remark}{Remark}

\newacronym{rfm}{RFM}{Ribosome Flow Model}
\newacronym{pmp}{PMP}{Pontryagin Maximum Principle}
\newacronym{tas}{TASEP}{Totally Asymmetric Simple Exclusion Process}

\newcommand{\mean}[1]{\overline{ #1}}

\newcommand*\diff{\mathop{}\!\mathrm{d}}
\begin{document}
	\maketitle 
	
\begin{abstract}
 A dynamical system \emph{entrains} to a periodic input if its state converges globally to an attractor with the same period. In particular, for a constant input the state converges to  a unique  equilibrium point for any initial condition. %
 We consider the problem of maximizing a weighted  average of the system's  output  along the periodic attractor.
 The \emph{gain of entrainment}  is the benefit  achieved by using a non-constant periodic input  relative to  a constant input 
 with the same time average.   Such a problem amounts to optimal allocation of resources in a periodic manner. We formulate this problem as a periodic optimal control problem which can be analyzed by means of the Pontryagin maximum principle or solved numerically via powerful software packages.
 We then apply our framework to a class of occupancy models that appear frequently in biological synthesis systems and other applications. We show that, perhaps surprisingly, constant inputs 
are optimal  for various architectures. This  suggests that the presence of non-constant periodic signals, which frequently appear in biological occupancy systems, is a signature of an underlying \emph{time-varying} objective functional being optimized.%
	\end{abstract}


	\section{Introduction}

	Periodic oscillations  are abundant  in biomolecular systems, and an extensive body of research has been devoted to study   their roles in intracellular and extracellular interactions \cite{fung05,tu06}.  In the presence of such excitations,
proper functioning of biological systems often  requires their internal states   to \emph{synchronize} with the periodic input signal.   In the parlance of systems theory, this is known as \emph{entrainment}, which means that the response of a system subject to a periodic input with period~$T$ will converge to a periodic trajectory of the same period~$T$. There has  been great recent interest in the study of this phenomenon \cite{ferrell11,entrainment,katz2020,entrainME,entrain_trans}. Examples of external periodic influences include  operation	under the influence of sunlight, which requires the internal 
	clocks of   biological organisms to entrain to the 24-hour solar day. For instance, 
	it has been shown that the plant  \emph{Arabidopsis}
	uses its circadian  clock 
	to anticipate times with an increased  susceptibility 
	to  fungal pathogens, 
	and regulates  its immune system resources
	accordingly~\cite{entrain_arabid}.  Entrainment is also essential in many \emph{synthetic} biological systems. For instance, synthetic oscillators can be used to emulate  natural hormone release rhythms in the treatment of certain diseases \cite{khalil10}. More generally,  robust and optimal     synthetic oscillators  constitute  an important  module in larger systems \cite{paulsson16,guan18}.
	
	At the intracellular level, the \emph{cell cycle} is a periodic routine that regulates DNA replication and cell division. This  requires precise regulation of many interacting proteins, and also  appropriate  resource allocation  at different stages of the cell cycle. Deviations from the program can lead to cell death or cancer. 
	
		An important underlying process   is translation, which is a major component in the central dogma of molecular biology, and   requires sophisticated coordination between  ribosomes, mRNA and tRNA molecules, and various proteins. Two of the key underlying steps are   \emph{initiation}
		in which the ribosome attaches to an mRNA molecule,
		and \emph{elongation},
	  in which the ribosome  scans along the mRNA  to produce a chain of amino-acids. 
		Regulation of initiation and elongation are an effective way to control protein concentrations \cite{sonenberg09,gobet17}. 
		One biological mechanism for  
		cell-cycle regulated genes is based on
		     codons whose corresponding tRNAs have low abundances (known as non-optimal codons)~\cite{frenkel12,zhou13}. In particular, \emph{periodic} variations in the level of these specific tRNAs can generate cell cycle-dependent oscillations in the corresponding protein levels \cite{frenkel12}. In other words,  the protein levels entrain to the periodic excitation provided by the tRNA levels.  Similar oscillation-inducing regulation mechanisms during DNA damage response have also been reported~\cite{patil12}. 
		Other works have indicated that the speed of translation is sensitive to fluctuating tRNA availability \cite{wohlgemuth13}, that cells use~tRNA to control protein abundance in stress conditions \cite{Torrenteaat6409}, and that tRNA disregulation is a contributing factor in cancer 
		progression~\cite{goodarzi16}. 
		In addition to tRNA regulation, many other intracellular oscillators have been identified  as regulators of the cell cycle \cite{cross03,da15}.
		
		In what follows, we  first describe, as a motivation,
		a  class of mathematical models that are useful in modeling various processes involved in gene translation. Our focus is to analyze  these models in the presence of periodic excitations modeled as periodic inputs. We then state the generic control problem to be solved.
		
		\subsection{Motivation: Occupancy models}
In many important biological models,   state variables  describe the  occupancy in a certain site or compartment. For example, in physiology
compartmental models describe  drug absorption distribution and elimination   in various
 body fluids or tissues~\cite{multi_comp2015}.

\subsubsection*{One-dimensional models}
	Many biological processes   involve ``biological machines'' that move	along a 1D lattice of ordered ``sites''.
Examples include ribosomes that scan mRNA during translation,
molecular motors that carry cargoes along
		a filamentous network in the cytoskeleton,
		and phosphotransferases that transfer the phosphoryl group from the sensor kinases to some  ultimate target. To be concrete, we focus on ribosomes and mRNA translation, but the same ideas apply to other models.

		We now derive such a 1D  occupancy model
		using  several alternative modeling approaches.
				Let~$X$ $(Z)$ be the species denoting bound (unbound) ribosomes, respectively. The free ribosomes bind to mRNA. Bound ribosomes need tRNAs to translate the information in the mRNA into proteins~$(P)$. 
				A phenomenological one-step model written in Chemical Reaction Network (CRN) formalism~\cite{erdi89} gives:
				\begin{align*} 
				{\rm mRNA } + Z & \to X \\
				{\rm tRNA} + X & \to Z + {\rm mRNA } + P. \vspace{-0.05in}
				\end{align*}
We assume that tRNA and mRNA are abundant, so that their dynamics are not affected by the above reactions. Note that the species~tRNA represents all possible variants of transfer~RNA.
				Let $x(t)$ be the concentration of occupied (bound) ribosomes in the cell at time $t$, and 
				let~$z(t)$ be the concentration of free ribosomes. The occupancy of ribosomes is determined by mRNA transcript abundance $u_0(t)$, and tRNA abundance $u_1(t)$.   The CRN gives the following system of bilinear ODEs:
				\begin{align} \label{2ode}
				\dot x(t) &= u_0(t) z(t) - u_1(t) x(t), \\ \nonumber
				\dot z(t) & = u_1(t) x(t) - u_0(t) z(t).
				\end{align}
				Assuming a fixed total concentration of ribosomes $M$, we have~$x(t)+z(t)\equiv M$. The total concentration can be normalized to $M=1$. Then the two-dimensional dynamics can be reduced to a one-dimensional ODE:
		\begin{equation}\label{1drfm_intro}\dot x(t) = u_0(t) (1-x(t)) - u_1(t) x(t).\end{equation}
				This implies that $x(t)$ evolves on the unit interval, and it can be interpreted as a normalized 
				occupancy of some site at time~$t$. 
				More generally,~$x(t)$
				can be interpreted as the  probability that a certain site is occupied  by some ``biological machine'' like a ribosome or a molecular motor.  This occupancy model has been termed a ``bottleneck'' module in \cite{RFM_IEEE_CL}.

				Note that the occupancy model~\eqref{2ode} can also be used to model binding and unbinding of a
				substrate to an enzyme.
				
		\subsubsection*{Multisite models: The Ribosome Flow Model}
		The \acrfull{tas}~\cite{solvers_guide} is a 
		fundamental stochastic  model from nonequilibrium statistical physics.  In~\acrshort{tas}, particles 
		move forward at random times along  a 
		1D chain of  sites. A site can be either   free or contain a single particle.
		Totally asymmetric means that the flow is unidirectional, and simple exclusion means that a particle can only hop into a free site. This models the fact that two particles cannot be in the same place at the same time.
		The simple exclusion paradigm generates  an indirect coupling between the particles, and also allows modeling the evolution of ``traffic jams'':  if a particle remains at site~$i$ for a long time, then particles will 
		accumulate ``behind'' it, i.e. in site~$i-1$, then site~$i-2$ and so on. TASEP has been used extensively to model and analyze ribosome flow~\cite{TASEP_tutorial_2011}
		and many more natural and artificial processes
		including molecular motors, 
		traffic flow, evacuation dynamics, and more~\cite{ScCN11}. 
		
		The \emph{Ribosome Flow Model} (RFM)~\cite{reuveni2011genome} is the dynamic mean-field approximation of~TASEP. In the~RFM,  the state-variables~$x_1(t),\dots,x_n(t)$ describe the occupancy in~$n$ sites along the mRNA molecule. The  RFM dynamics is described by a system of~$n$ first-order ODEs:
	\begin{align}\label{eq:rfm}
		\dot x_k &= \lambda_{k-1} x_{k-1}  (1-x_k ) 
		-  \lambda_k x_k  (1-x_{k+1} ), \; k=1,\dots,n,
		\end{align} 
		where we define~$x_0(t)\equiv 1$ and~$x_{n+1}(t)\equiv 0 $. Here~$x_i(t)$ describes the occupancy at site~$i$ at time~$t$, normalized such that~$x_i(t)=0$ [$x_i(t)=1$] means that site~$i$ is completely 
		empty [full] at time~$t$.  
		In the context of   translation,~$\lambda_i(t)>0 $  describes
		the transition rate from site~$i$ to site~$i+1$ at time~$t$. This rate depends 
	  on various biomechanical properties, for example,
		the abundance of tRNA molecules  delivering the amino-acids to the ribosomes.
		Eq.~\eqref{eq:rfm} can be explained as follows. The change in the density in site~$k$ is the 
		flow from site~$k-1$ into site~$k$ minus the flow from site~$k$ to site~$k+1$.
		The first term,~$\lambda_{k-1} x_{k-1}  (1-x_k ) $,   is proportional to the transition rate from site~$k-1$ to~$k$, the occupancy at site~$k-1$, and the amount of ``free space''~$(1-x_k)$ at site~$k$. 
		Note that this is a ``soft'' version of   simple exclusion. 
		The second term is similar. 
	Note that~$\lambda_n(t) x_n(t)$ describes the flow of ribosomes out of the last site at time~$t$, i.e. the protein production rate. If the whole mRNA strand is considered as one site, that is,~$n=1$ then the RFM model will be identical with the occupancy model~\eqref{1drfm_intro}.

		The state-space of the~RFM is the~$n$-dimensional  unit cube $[0,1]^n$. It was shown in~\cite{margaliot2012stability} 
		that the RFM (with constant $\lambda_i$'s) admits a unique equilibrium~$x_e=x_e(\lambda_0,\dots,\lambda_n)
		\in(0,1)^n$, and that for any~$a \in[0,1]^n$ the corresponding solution
		of~\eqref{eq:rfm} satisfies~$\lim_{t\to \infty} x(t;a)=x_e$. In other words, the transition rates
		determine a unique Globally Asymptotically Stable (GAS) equilibrium.  	More generally, Ref.~\cite{entrainment} showed that if the rates
		are time-varying,  and jointly periodic with a period~$T$, then~\eqref{eq:rfm} admits a 
		GAS solution~$\gamma_T :\R_+ \to   (0,1)^n$, that is~$T$-periodic,
		and~$x(t,a)$ converges to~$\gamma_T$ for all~$a\in   [0,1]^n$. In other words, the RFM entrains. 
		Note that a constant rate is $T$-periodic for any~$T$, so
		entrainment also holds if   a single rate   
		is~$T$-periodic and all the other rates are constant.
		In the biological  context, entrainment can be interpreted as follows: if, say, variations in tRNA abundances generate  $T$-periodic initiation and/or elongation rates, then the protein production rate will also converge to a periodic pattern with period~$T$.

		The RFM and its variants have been used extensively to model and analyze
		ribosome flow during the process of translation (see e.g.~\cite{nani,alexander2017,Raveh2016,rfm_max}),
		as well as other important cellular processes like phosphorelay~\cite{EYAL_RFMD1}.

		Just like~\acrshort{tas},
		 the RFM (and in particular the model~\eqref{2ode})
		is a phenomenological model that can be applied to 
		study various processes like vehicular or pedestrian traffic~\cite{RFM_IEEE_CL}. In this case, the occupancy is interpreted as the ratio between  the number of vehicles (or pedestrians)    at a certain junction at time~$t$ and the total number of possible vehicles.

		\subsubsection*{Generalized occupancy models}
		Let~$\R^n_+:=\{x\in \R^n \st x_i \geq 0,\; i=1,\dots,n\}$ denote the non-negative orthant in~$\R^n$. 
		Recall that the linear single-input single-output~(SISO) linear system
		\begin{align*}
		\dot x&=A x +bu, \\
		y&=c^Tx,
		\end{align*}
		is called \emph{positive} if every entry of~$b,c$,
		and every off-diagonal entry of~$A$ is non-negative (that is,~$A$ is a \emph{Metzler} matrix).
		This implies that for any~$x(0)\geq 0$ and any control~$u$ with~$u(t)\geq 0$ for all~$t$ we have~$x(t)\geq 0$ and~$y(t)\geq 0$ for all~$t\geq 0$~\cite{farina2000}. This is useful when the state-variables and the output represent physical quantities that can never attain negative values, e.g., population sizes or concentrations of molecules.

		Generalized occupancy models~(GOMs) are   a cascade of occupancy models and   SISO positive  linear systems. 
		These models are useful when the output of an occupancy model is the input to another biological system, that in the vicinity of its equilibrium point, can be approximated as a positive linear dynamical system. 
		Similar to the multisite RFM model introduced before, it can be shown that GOMs entrain to periodic inputs.	
			
\begin{figure}
			\includegraphics[width=\textwidth]{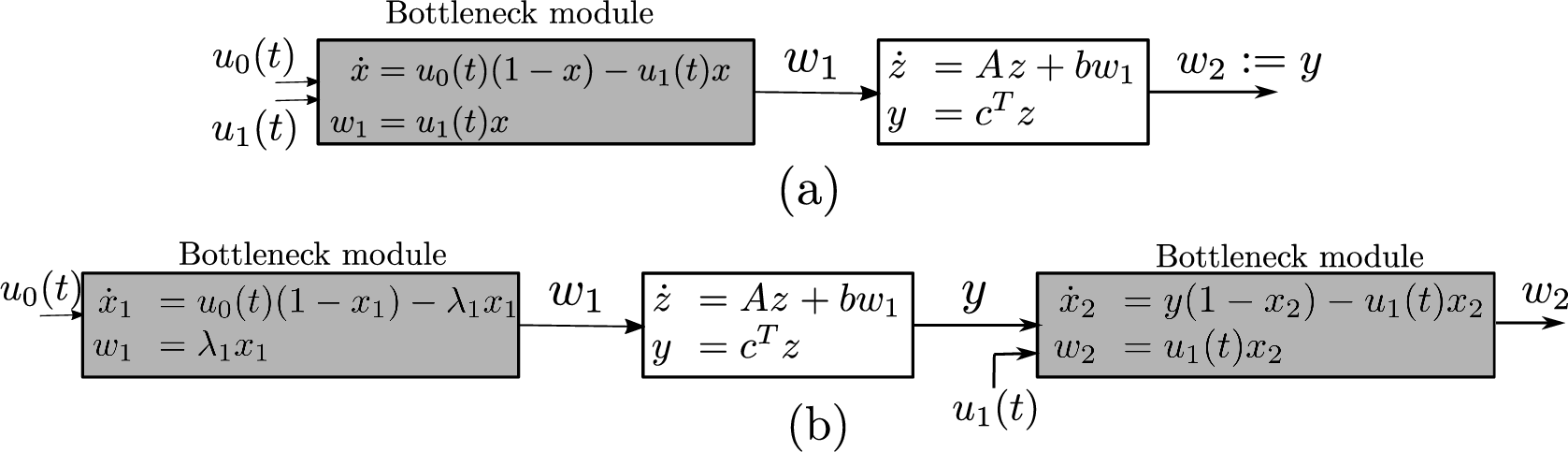}
			\caption{Two examples of generalized occupancy models . The controls are $u_0(t),u_1(t)$ which are scalar functions. We have $x_1, x_2, x, w_1, w_2, y \in \mathbb R_{+}$, $z \in \mathbb R^n, A \in \mathbb R^{n \times n}, b,c \in \mathbb R_{+}^n$. The linear system block is assumed to be positive and Hurwtiz.   }
			\label{f.cascade}
		\end{figure}

			For example, 		Fig.~\ref{f.cascade}-(a) depicts a time-varying bottleneck module feeding a positive linear system. In this module,~$u_0,u_1$ are entrance rates, and~$w_1$ is the exit rate. The effective inflow is proportional to the \emph{vacancy}~$1-x(t)$, while the outflow is proportional to the \emph{occupancy}~$x(t)$.  This  cascade  models an occupancy model driving a downstream linear system. 
	As another example, 	  Fig.~\ref{f.cascade}-(b) depicts a linear system ``sandwiched'' between a 2-site RFM and 1-site RFM. This can model a situation where the
production rate of one protein affects, via another biological process, the promoter (and thus the transcription initiation rate) of some other  mRNA.

		A GOM can also be used to model the RFM with  time-varying rates 
		under the condition
		\be\label{eq:adftp}
		\lambda_i(t) \gg \lambda_0(t) \text{ for all } i\geq 1 \text{ and all }t\geq 0.
		\ee
		Then we can expect that the initiation rate becomes the bottleneck rate and thus~$x_i(t)$, $i=2,\dots,n$, converge to values that are close to zero, 
		suggesting that~\eqref{eq:rfm} can be simplified to 
		\begin{align}\label{eq:rfm_l}
		\begin{split}
		\dot x_1(t) &= \lambda_0(t)  (1-x_1(t) )  - \lambda_1 (t) x_1 (t) , \\
		\dot x_i (t) &= \lambda_{i-1} (t) x_{i-1} (t) -  \lambda_i(t) x_{i},\;\; i\in \{2,\dots,n\}, 
		\end{split}
		\end{align} 
		which has the same form as the cascade in Fig.~\ref{f.cascade}-(a).

After these motivating examples, we next formulate the abstract questions to be studied in this paper. 
		
		\subsection{Gain of entrainment}

		Entrainment  can be studied  in the framework of systems and control theory. The periodic excitation is modeled as the control input~$u(t)$ of  a dynamical system,  and the system entrains if in response to a~$T$-periodic excitation 
		it    admits a globally attractive~$T$-periodic solution~$\gamma_T$.
		In other words, every solution of the system converges to the attractor~$\gamma_T$.

		Here, we consider a   \emph{quantitative}  potential 
		advantage of entrainment 
		called  the \emph{gain of entrainment}. To explain this, consider a control system 
		that, for any~$T\geq 0$ 
		and any  $T$-periodic control~$u_T$,
		entrains to a  unique~$T$-periodic solution~$\gamma_T$.  
		Note that, in particular, 
		this implies that for any constant control~$u(t)\equiv u_0$
		the trajectory converges to a unique equilibrium~$\gamma_0$ for any initial condition. 
		Suppose also that the system admits a  scalar output~$y(t)=h(t,x(t),u(t))$ (that is, a function of time, state, and  the input), 
		and that~$h$ is~$T$-periodic in~$t$, so that the output also entrains. 
		The output represents a quantity that we would like to maximize, e.g.~traffic flow or protein production rate.
		
		Since the system entrains, we ignore the transients and consider the problem of maximizing the average of the periodic output, that is, the average over a period of~$h(t,\gamma_T,u_T)$.
		The gain of entrainment is the benefit (if any) 
		in the maximization for a (non-trivial) periodic control
		over a constant control.  A natural example is  to analyze
		the gain in  traffic flow  for periodically-varying traffic signals over
		constant signals. However, to make this meaningful, we must add another assumption, namely,  that the total 
		time of green lights in both alternatives is equal. 
		Mathematically, this means that we compare the average output for a time-periodic control~$u_T$
		and a constant control~$\bu$ such that the average value
		of~$u_T$ over a period is equal to $\bu$. 
		If  the gain of entrainment is positive then entrainment does not only assist in producing   an internal  clock that can follow an external periodic excitation, but also 
		yields higher production rates than those obtained by equivalent constant excitations.

		The possible advantages of
		periodic forcing of various production processes are well-known.
		For example, Ref.~\cite{kol_perio} states that:
		``...theoretical and experimental studies
		have shown that the performance (for instance micro-algae or bio-gas production) of some optimal steady-state
		continuous bioreactors can be improved by a periodic modulation of an input such as dilution rate or
		air flow''.	Ref.~\cite{per_fish} studies a PDE model for harvesting a biological resource  and  demonstrates the advantages of periodic harvesting over  a constant one.

			 The   gain of entrainment was recently introduced in~\cite{RFM_IEEE_CL}. 
		Entrainment in nonlinear systems is nontrivial to prove. 
		A  typical proof is based on  contraction theory~\cite{entrain_trans,entrainment}, yet this type of proof
		provides no information on the attractive periodic solution,  except for its period  (see~\cite{coogan_margaliot} for some related considerations). 
			 Nevertheless, we show here that  determining the  gain of entrainment can be cast as an \emph{optimal control problem}.
This allows using powerful theoretical tools, like Pontrayagin's maximum principle~\cite{pontryagin,LeeMarkus,liberzon}, as well as  numerical methods  in studying the gain of entrainment. We demonstrate this by analyzing the gain of entrainment in several examples of  occupancy models. 
	 
	  For instance, consider the gain of entertainment for~\eqref{eq:rfm_l}.  It is natural to speculate that using time-periodic 
		rates~$\lambda_i(t)$, that are properly synchronized,   yields a positive gain of entertainment with respect to using constant rates (with the same average values). 
		In the context of traffic flow, this is equivalent to the conjecture that   properly synchronized 
		periodic traffic lights  can improve the overall flow. However, we  show that, perhaps surprisingly, for a subclass of these systems the  gain of entrainment is zero.  
		
		We also consider a  problem formalism that allows for
		time-varying costs of resources, like tRNAs, along the period. These may be produced at different unit costs at different times of the cycle. This modified formulation allows the allocation of resources differently at different times along  the cycle. Also, instead of average throughput,   a \emph{weighted} average of the product may be more relevant, in the sense that we may need certain enzymes at different times of the day or at different points in the cell cycle.  This 
		corresponds to 
		``just-in-time production''~\cite{zaslaver04}.
		In such cases, we show, not so surprisingly, 
		that time-varying periodic inputs may indeed offer an advantage over
		 constant inputs.
		This  suggests that the presence of non-constant periodic signals, which frequently appear in biological occupancy systems, implies that the system is optimising 
		an underlying \emph{time-varying} objective functional.
		
		Our work is   related to results from the  field of optimal periodic control~(OPC) (see, e.g.,~\cite{colo88}). 
		As noted by Gilbert~\cite{Gilbert77},  OPC was 
		motivated by the following question:
		Does time-dependent periodic control yield better process performance than optimal
		steady-state control? 
		In particular, the recent paper~\cite{over_yield} defines a notion called \emph{over-yielding} that is 
		closely related to the gain of entrainment. 
		However, our
		setting is different, as in~OPC periodicity was enforced by restricting attention to controls~$u$ guaranteeing that~$x(T)=x(0)$. This implies in particular that the initial value~$x(0)$ (and thus also general transient behaviors)
		may have a strong effect on the results.
 Also,  in the typical
  OPC formulation
		there is in general 
		no requirement that the averages of the periodic and constant controls are equal.

		We study systems that entrain and thus
		for a~$T$-periodic control  the state of the system  converges to a unique~$T$-periodic trajectory for \emph{any} initial condition~$x(0)$. In other words,  we  consider  the behavior of attractors.

  The remainder of this paper is organized as follows. The next section defines
	the gain of entrainment for a general 
 mathematical model. Section~\ref{sec:otcon} shows 
how the analysis can be cast as an optimal control problem. 
 Section~\ref{sec:1drfm} demonstrates the theory for the two-input bottleneck module.  Section~\ref{sec:entrain_GOM}
 proves that for several GOMs, including the ones depicted in Fig.~\ref{f.cascade},
the gain of entrainment is zero. Finally, conclusions and future directions are presented in Section~\ref{sec:conc}.
The appendix contains proofs of the results including a detailed analysis characterizing extremals via the Pontraygin Maximum Principle (PMP).

 \section{Gain of entrainment }
We 	consider a general nonlinear  control  system:
	\begin{align}\label{eq:dyns}
	\dot x&=f(x,u),\\
	y&=h(t,x,u),\nonumber 
	\end{align}
	with $f,h$  locally Lipschitz functions, the state~$x(t)\in\R^n$, control (or input) $u(t)\in\R^m$,
	and scalar output~$y(t)\in\R$. We allow $h$ to be time-varying to include the cases in which  different weights can be used at different times in the cycle.
	The set of admissible controls consists of
	measurable functions taking values in some closed and compact set~$U\subset \R^m$. 
	Let~$x(t,p,u)$ denote the solution of~\eqref{eq:dyns} at time~$t\geq 0$ for the initial condition~$x(0)=p$ and the control~$u$. 
	We assume throughout that 
	for any~$x(0)$ in the state-space and any admissible control, \eqref{eq:dyns} admits a unique solution  for all~$t\geq 0$.

	We say that  system~\eqref{eq:dyns}   \emph{entrains}
	if in response to 
	any admissible and~$T$-periodic control~$u_T$ the system  admits a unique $T$-periodic solution~$\gamma_T(t)$ (that depends 
	on~$u_T$), and
 for any initial condition~$p$ the solution~$x(t,p,u_T)$ converges to~$\gamma_T$. 
This  implies in particular that
the system ``forgets'' its initial condition.  
To explain the mathematical formulation of the gain of entrainment, fix~$q\in\R^m$ with~$q_i>0$ for all~$i$. 
We would like to consider only 
inputs whose average over a period is~$q$, and compare their effect to the effect of the constant control~$u(t)\equiv q$. However, we allow a slightly more general scenario by fixing  a weighting function~$\alpha(t)>0$ such that $\frac{1}{T}\int_{0 } ^T \alpha(t) \, \diff t =1$. We then restrict
attention 
 to~$T$-periodic controls 
satisfying the weighted integral constraint:
\be\label{eq:intcons}
				\frac{1}{T} \int_{0}^T \alpha(t) u (t)  \diff   t =q, 
\ee
that is, the $\alpha$-weighted average   of~$u $ is~$q$. This can be further 
generalized by allowing a general measure $\mu$ on the interval $[0,T]$ and imposing $\int_{[0,T]} u (t) \diff \mu = q$. However, we keep the presentation simple by adhering to \eqref{eq:intcons}.

Let
\be\label{eq:zut}
z(u):=\frac{1}{T} \int_{0}^T h(t,\gamma_T(t),u (t)  )  \diff  t , 
\ee
that is, the   average  
value of the output along the globally attractive~$T$-periodic solution (recall that we assume that~$h$ is~$T$-periodic in its first variable). 
If the convergence to~$\gamma_T$ is relatively fast then
after a short transient the average output over a period of length $T$ is very close to~$z(u )$.
In applications in fields like biotechnology   and traffic control  
the average value of the output, and not its specific
values at all times,  is often the relevant quantity.

The constant control~$u(t)\equiv q$, which we simply 
denote by~$q$, is also~$T$-periodic (for any~$T\geq 0$) and satisfies~\eqref{eq:intcons}. Hence, the corresponding solution converges to a fixed point~$e=e(q)$
and 
\begin{align*}
z(q) &=\frac{1}{T} \int_{0}^T h(t,e,q )  \diff   t\\&=h(e,q) . 
\end{align*}
The \emph{gain  of entrainment} of~\eqref{eq:dyns} is  defined as
\begin{equation}\label{gain}
				c_T(q):=    \sup_{u } z(u )    -z(q), 
\end{equation}
where the sup is over all admissible, $T$-periodic controls
 that satisfy the constraint~\eqref{eq:intcons}. Thus, we are always
comparing the effect of controls with the same average value.  Note that~$c_T(q) \ge 0$ for all $q$.
If~$c_T(q)>0$ for some $q$, then there exists a nontrivial periodic control that yields
a higher average output than that obtained for  a constant control.
If~$c_T(q)=0$, then nontrivial~$T$-periodic controls are ``no better''
 than the simple constant control equal to $q$.

	To gain a wider perspective, consider the case of 
	a~SISO asymptotically stable~LTI system
	with input [output] $u(t)$ [$y(t)$]
	and transfer function~$G(s)$. 
	Fix~$T>0$, and
	consider the ~$T$-periodic control
	\[
	  u_T(t):=  a+b\sin(2\pi t/T),
	\]
	with~$a,b\in \R$.  
	Note   that~$\frac{1}{T}\int_0^T u_T(t) \diff  t =a$. Let~$\omega:=2\pi/T$.
	It is well-known that   the output converges 
	to~the~$T $-periodic function~$y_T(t):= G(0)   a+
	|G(j\omega)| b   \sin( \omega t + \angle G(j \omega))$, where~$j:=\sqrt{-1}$,
	so~$  \frac{1}{T}\int_0^T y_T(t) \diff   t   = G(0)  a  $. 
	On the other-hand, for the constant control~$u(t)\equiv a$  
	  the output converges to~$ G(0)  a$, which is the same value. Thus, for this input the gain of entrainment is   zero.
		Any~$T$-periodic, measurable, and bounded input  can be expressed as a Fourier series in terms of sinusoidal functions,  and this implies 
		that for LTI  systems the gain of entrainment is always zero.

\begin{figure}[t!]
		\centering
		\includegraphics[width=.6\linewidth]{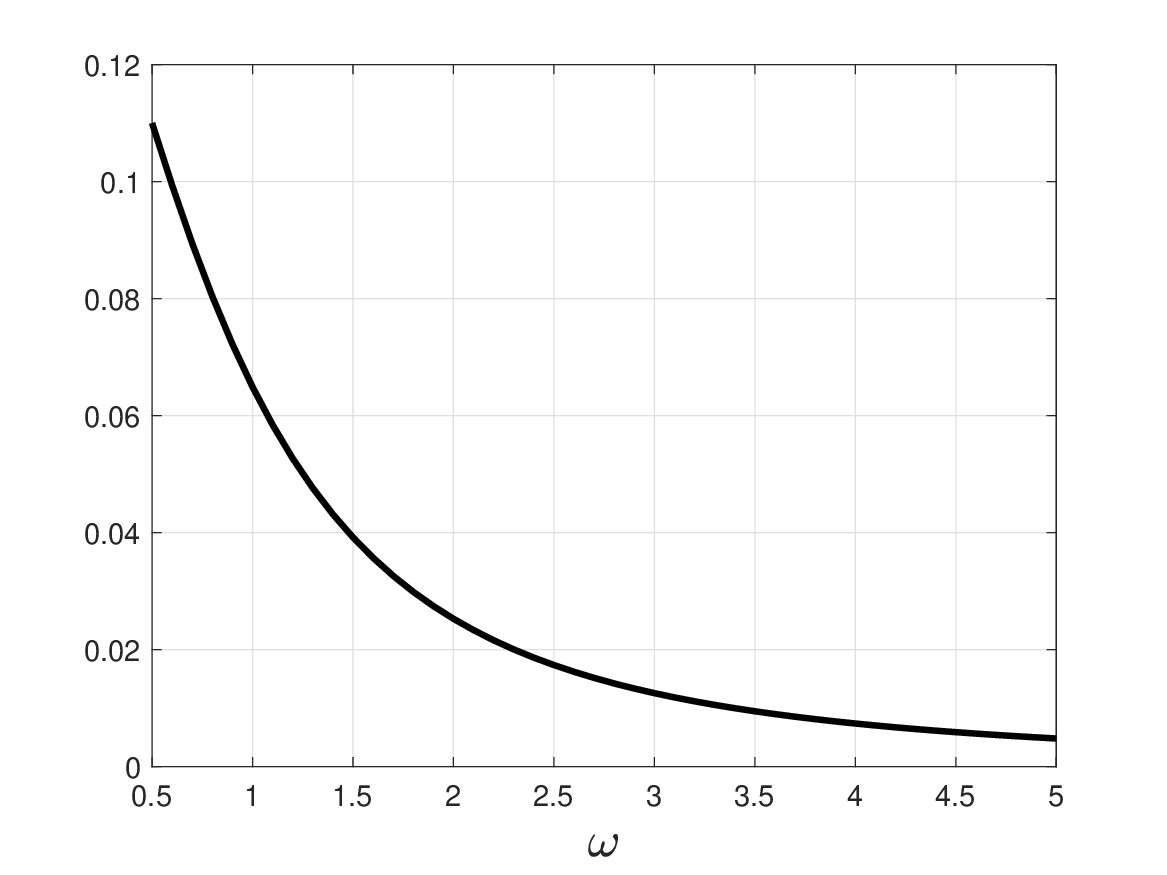}
		\caption{Average value of~$x_T(t)-1$ in Example~\ref{exa:sinsimp}
		as a function of~$\omega \in[0.5,5] $.}
		\label{fig:sin}
	\end{figure}

	However, for nonlinear system 
	the   gain of entrainment 
	may be positive. 
	 The next  two examples demonstrate  this.
	\begin{Example}\label{exa:sinsimp}
		Consider the scalar system:
		\begin{align}\label{eq:smsi}
		\dot x(t)&= 1-x(t) u(t), \\
		y(t)&=x(t).\nonumber
		\end{align}
		Fix~$T>0$. For a function~$v:\R_+ \to \R_+$, let~$\overline  v:=\frac{1}{T} \int_0^T v(s)\diff s$.
		Fix~$q>0$. For the control~$u(t)\equiv q$ any  solution of~\eqref{eq:smsi}  converges  
		to the equilibrium~$q^{-1}$. Consider a~$T$-periodic and positive control~$u_T(t)$ satisfying~$\bar {u}_T=q$, and assume there exists some $\alpha>0$ such that $u_T(t)\ge \alpha$ for almost all $t \in[0,T]$. Then any matrix measure  of the Jacobian of~\eqref{eq:smsi} is uniformly less or equal than $-\alpha<0$. Therefore, the system is contractive and any solution of~\eqref{eq:smsi}  converges to a unique~$T$-periodic solution~$x_T(t)$.
Let~$\omega:=2\pi/T$. 
		Consider now  the specific~$T$-periodic  control
		\be \label{eq:perub}
		u_T(t):=1 +(1/2) \cos(\omega t).
		\ee
		Here, $q=\tfrac 1T \int_0^T u_T(t) \diff t=1$. For this input, 
		the corresponding solution of~\eqref{eq:smsi} is:
		\[
		x(t)=\exp\left(-t-\frac{\sin(\omega t)}{2\omega } \right) (x(0)+ \phi(t)) ,
		\]
		where~$\phi(t):=\int_0^t \exp\left (s+\frac{ \sin(\omega s)}{2\omega} \right ) \diff s$. In particular,  
		\[
		x(T)=\exp(-T  )
		(x(0)+ \phi(T) ).
		\] 
		The initial condition~$x(0)=c$ for which the solution
		is~$T $-periodic  is
		\[
		c=\exp(-T) (c+ \phi(T)),
		\]
		so
		\begin{equation}\label{c_example}
		c=\frac{\exp(-T) \phi(T)}{ 1-\exp(-T)}.
		\end{equation}
		Thus, the attractive periodic solution is  
		$x_T(t):= \exp\left(-t-\frac{ \sin( \omega t)}{2\omega} \right ) (c+ \phi(t)) $.

		The average of the control~$u_T(t)$ is $q=1$.
		On the other-hand, for  the control~$u_0(t)\equiv 1$
		the solution of~\eqref{eq:smsi}
		converges to the steady-state~$1$.
		Fig.~\ref{fig:sin}  depicts the value
		\[
		    \frac{1}{T} \int_0^{T } x_T(t) \diff t  -1,
		\]
	 as a function of~$\omega=2 \pi /T$. It may be seen that this is always positive, and is maximal as~$T \to \infty$.

	We conclude that for~$q=1$ the
			gain of entrainment of~\eqref{eq:smsi} is  positive
			for any~$T >0$. Note that for large values of~$\omega $
			the gain of entrainment goes to zero. This is expected due to \emph{averaging}~\cite{khalil_third_edition_2002}. Roughly speaking, for large values of~$\omega$ the system cannot track the fast changes in the input, and thus responds to the average of the input.
			More rigorously,  for a system 
			   affine in the control, the map from controls on an interval~$[0,T]$ to trajectories on~$[0,T]$ is continuous with respect to the weak$^*$ topology in~$L^1$ and the uniform topology on continuous functions, respectively (see, e.g.~\cite[Theorem~1]{sontag_book}), and (2) for a periodic input~$u(t)$, the input~$u(\omega t)$ converges weakly to the average of~$u$.  An alternative proof is given for example in the textbook~\cite{khalil_third_edition_2002} (Section 10.2) (changing time scale in the statement of Theorem~10.4, by $x(t) = x(t/\epsilon)$). 
	\end{Example} 
	 
	\begin{Example} Consider the   system:
	\begin{align} \label{eq:consta}
	\dot x_1(t)&= - x_1(t) + u(t),\nonumber \\
	\dot x_2(t) &= -x_2(t) +  a x_1^2(t),\\
	y(t)&=x_2(t),\nonumber
	\end{align}
	with~$a>0$.
	
	Consider the input~$u_T(t):=1+\sin(\omega t)$, with~$\omega >  0$. Here, $q=\bu=1$. Let~$T:=2\pi/\omega$, and let~$H(s):=\frac{1}{1+s}$.
	Then~$x_1$ converges to the steady-state solution:
	\[
	 x_{1T}(t) := 1+|H(j\omega)| \sin(\omega t  +\angle H(j\omega) ).
	\]
     Hence,
	\begin{align*}
	 x_{1T}^2(t)&=
	1+  2|H(j\omega)| 
	\sin (\omega t  +\angle H(j\omega)  )+  |H(j\omega)|^2 \sin^2(\omega t   +\angle H(j\omega) )\\
	&=	1+ \frac{|H(j\omega)|^2 }{2 }+   2|H(j\omega)| \sin (\omega t  +\angle H(j\omega))
	-\frac{|H(j\omega)|^2  }{2 } \cos(2\omega t + 2 \angle H(j\omega)).
	\end{align*}
	It   follows that~$x_2$ converges to the steady-state solution
	\begin{align*}
	 x_{2T}(t) :=& a\left(
		1+ \frac{|H(j\omega)|^2  }{2 } +
	  2  |H(j\omega)|^2  \sin (\omega t  + 2\angle H(j\omega) ) \right.\\
	  & \left. -\frac  {|H(j\omega)|^2  }{2}|H(2\omega)| \cos(2\omega t  + 2\angle H(jw) +
	  \angle H(2j\omega ) )\right ),
	\end{align*}
	so
	\[ \bar x_{2T}=
	\frac 1T \int_0^T x_{2T}(t) \diff t =a\left(
		1+ \frac{|H(j\omega)|^2  }{2 } \right).
	\]
	 On the other-hand, for the average input~$\bar u=\frac{1}{T}\int_0^T (1+\sin(\omega t))\diff t =1$,
	 $x_1(t)$ converges to one, and
	 $x_2(t)$ to~$a$, so the average of the output is~$a$.
	 The difference between the two averaged outputs is thus
	\[
	\frac{a}{2}|H(j\omega)|^2=
	\frac{a}{2(1+\omega^2)} .
	\]
	This is maximized for~$\omega=0$, so the gain of entrainment is at least $c_T(1)=a/2$. 
	Observe that examples with arbitrarily large gain of entrainment can be obtained by 
	 taking  the constant~$a$ in~\eqref{eq:consta} large enough.
\end{Example}

In the next  section, we cast the problem of determining the gain of entrainment as an optimal control problem. 

	\section{Optimal control formulation} \label{sec:otcon}
Consider the control 
	system~\eqref{eq:dyns} with~$n$ state-variables and~$m$ inputs.
	We assume that the system entrains. 
	Pick any~$T>0$ and any~$q \in \R^m$. 
	We restrict attention to~$T$-periodic controls
	satisfying the individual weighted average constraints:
	\be\label{eq:uxint}
	\frac{1}{T}\int_0^T \Xi (t,u (t)) \diff t=q,
	\ee
	where $\Xi:\mathbb R^m_{+} \times [0,T] \to \mathbb R^m_+ $ is an integrable vector
	positive function that satisfies $\frac{1}{T}\int_0^T \Xi(t,q) \diff t=q$.

The $n$-dimensional control system with the integral constraint on the  controls \eqref{eq:uxint} can be lifted to
an $(n+m)$-dimensional nonlinear control system by adding $m$-equations to \eqref{eq:dyns}:
\begin{equation}\label{lifted_sys}
 \dot{\tilde x}=\begin{bmatrix} \dot x \\ \dot{\xi}   \end{bmatrix}  = \begin{bmatrix} f(x,u) \\ \Xi(t,u) \end{bmatrix}=F(t,\tilde x,u ),
\end{equation}
where $\tilde x:=\begin{bmatrix}x^T  & \xi^T\end{bmatrix} ^T$. We impose the
boundary conditions:
\begin{align}\label{boundary_general}
x(0)=x(T), \; \xi(0)=0, \; \xi(T)=Tq .
\end{align}
Since we consider systems that entrain, for any~$T$-periodic control
there corresponds  a unique  GAS $T$-periodic solution~$\gamma_T(t)$. 
	The condition $x(0)=x(T)$ guarantees that the maximization is performed over this solution.
  The other two conditions are equivalent to~\eqref{eq:uxint}.

To make the problem well-posed, we will assume that   controls take values in a hypercube~$[\ell,L]^m$ where~$0< \ell<L$.
We then formulate an optimal control problem as follows:

\begin{problem} \label{prob:optim0}
	Fix values $0<\ell < q < L$. Find an admissible control~$u$ that maximizes the objective functional 
	\[
	J(u ):=\frac{1}{T} \int_{0}^T h( t, x(t),u (t) )  \diff  t , 
	\]
	subject to the ODE~\eqref{lifted_sys}, the constraint~\eqref{eq:uxint}, $u_j (t)\in [\ell,L], j=1,\dots,m, t \in[0,T]$, 
	and the boundary 
	conditions~\eqref{boundary_general}.
\end{problem}
Note that~$J(u)$ is    the   average  
value of the output along the globally attractive~$T$-periodic solution (recall that we assume that~$h$ is~$T$-periodic in its first variable). 

In what follows we always consider systems affine in the control. Then the fact that $[\ell,L]^m$ is compact and convex implies, by Filippov's Theorem (see, e.g.~\cite{agrachev-book}), that the reachable set at any time~$t\geq 0$ is compact. Since~$h$ is locally Lipschitz, an optimal control exists. 

The optimal control  formulation allows  
to apply powerful theoretical tools for solving optimal control problems as well as utilize 
software  packages for numerical solutions~(e.g.,~\cite{gpops,falcon_man}). 
In  the next section,
we   demonstrate how to determine the gain of entrainment using
this   formulation  for both time-invariant and  time-varying cost functions.

\begin{remark}
	The control signals can be assumed to belong to different intervals. In other words, we can have $u_j(t) \in [\ell_j,L_j]$, $j=1,\dots,m.$ Nevertheless, we have simplified the formulation above by re-scaling the controls so  that they all satisfy  the same bounds.  
\end{remark}

The following result is immediate:
\begin{theorem}\label{recast_as_opc}
    A control $u$ is a solution of Problem \ref{prob:optim0} iff it maximizes $c_T(u)=z_T(u)-q$ as defined in \eqref{gain}.
    \end{theorem}
    In other words, to find the gain of entrainment we must find a control~$u$  that solves  Problem~\ref{prob:optim0} and then compute~$J(u) - q$.

\subsection{Pontryagin's Maximum Principle (PMP)}
Problem \ref{prob:optim0} can be studied in the framework of the PMP~\cite{pontryagin,LeeMarkus,liberzon}. The \textit{Hamiltonian} associated with our problem is:
\begin{align} \label{eq:hamiltonian_general}
\mathcal H(t,u,\tilde x,p,p_0 ) &:= p^T(t)F(t,\tilde x,u )
+  \tfrac{p_0}{T} h(t,x,u ) ,
\end{align}
where $ p(t) \in \mathbb R^{n+m}$ is the co-state, and 
the \textit{abnormal multiplier}
$p_0\ge 0$ is a constant. 

\begin{proposition}[PMP] \label{p.pmp}
	Let~$u^*(t) \in \mathbb R^m , t \ge 0$ be an optimal control for 
	Problem~\ref{prob:optim0}, and
	let~$\tilde x^*:[0,T] \to \R^{n+m} $ be the corresponding optimal   trajectory.  There exist $p_0^* \ge 0$ and~$p^*: [0,T] \to  \R^{n+m} \setminus\{0\}$,  such that:
	\begin{enumerate}
		\item The optimal state~${\tilde x}^*(t)$ and corresponding adjoint $p^*(t)$ satisfy:
		\begin{align} \label{pdot}
		\dot  {\tilde x}^* & = \frac{\partial \mathcal H}{\partial p}(t,u^*,\tilde x^*,p^*,p_0^*) ,\nonumber \\ 
		\dot  p^* & = -\frac{\partial \mathcal H}{\partial \tilde x}(t,u^*,\tilde x^*,p^*,p_0^*).
		\end{align}
		\item The control $u^*(t)$  satisfies
		\begin{equation}\label{eq:pmp_max}
		\mathcal H(t,s,\tilde x^*(t),p^*(t),p_0^* ) \le \mathcal H(t,u^*(t),\tilde x^*(t),p^*(t),p_0^* )
		\end{equation}
		for all $s \in [\ell,L]^m$  and almost every (a.e.)~$t \in [0,T]$ .
		\item The adjoint  satisfies the 
		transversality condition:
		\begin{equation}\label{eq:trans}
		p_i^*(0) = p_i^*(T), \; i=1,..,n.
		\end{equation}
		\item $H(t,u^*(t),\tilde x^*(t),p^*(t),p_0^* )=0$  for all $t \in [0, T ]$.
	\end{enumerate}
\end{proposition}

\noindent		The proof is given in the Appendix. Utilization of the PMP to deduce the structure of the optimal control is a difficult problem in general. We will show in the next section how it can be utilized in certain cases.

	A trajectory $\mathscr X:=(u(t) ,\tilde x(t), p(t))$ is said to be \emph{feasible} if it satisfies the ODEs~\eqref{pdot} and the boundary conditions~\eqref{boundary_general} and~\eqref{eq:trans}. A feasible trajectory~$\mathscr X$ is an \textit{extremal trajectory} if it satisfies the~PMP, i.e.\ if it also satisfies Proposition \ref{p.pmp}. Observe that any optimal trajectory must be an extremal by Proposition~\ref{p.pmp}. %
	 \section{Occupancy models with controlled inflow  and outflow }\label{sec:1drfm}
	In this section we look more closely at the $n$-dimensional occupancy model of the form:
	\begin{align}  \label{rfm1}
 \begin{bmatrix} \dot x_1   \\ \dot x_2 \\ \vdots \\ \dot x_{n-1}\\ \dot x_n \end{bmatrix} = \begin{bmatrix} u_0(t) (1-x_1 ) - \lambda_1 x_1 (1-x_2)  \\ \lambda_1 x_1 (1-x_2) -\lambda_2 x_2 (1-x_3) \\ \vdots  \\  \lambda_{n-2} x_{n-2} (1-x_{n-1} ) -  \lambda_{n-1} x_{n-1} (1-x_n )  \\  \lambda_{n-1} x_{n-1} (1-x_n ) - u_1(t) x_n \end{bmatrix}.
\end{align}
This is an~$n$-dimensional RFM with 
    initiation  and exit rates that are non-negative control inputs. 	Suppose that both~$u_0(t),u_1(t)$ are periodic with period~$T\geq 0$. It was proved in~\cite{entrainment} that the RFM with $T$-periodic rates entrains, 
    so in particular~\eqref{rfm1} 
	admits   a unique solution~$\gamma_T(t)$, with~$\gamma_T(0)=\gamma_T(T)$,
	and~$x \to \gamma_T$ for any initial 
	condition~$x(0)\in[0,1]^n$.

	To study the cost of entrainment in this system,  
	fix~$0<\ell<L$. 
	We will assume that the rates~$u_0(t)$, $u_1(t)  $ are two measurable and essentially locally bounded functions
	taking values in the interval~$[\ell,L]$.
	
	To allow a ``fair''  comparison between~$T$-periodic controls and constant controls,  fix two values $ \bu_0, \bu_1 \in (\ell,L)$, and pose   integral constraints on the controls: 
\be\label{eq:oindd}
\frac{1}{T} \int_0^T \alpha_0(t) u_0 (t) \diff t=  \bu_0, \ 
\frac{1}{T} \int_0^T \alpha_1(t) u_1 (t) \diff t= \bu_1, 
\ee
for some given positive measurable functions $ \alpha_0(t), \alpha_1(t)$ satisfying   
\[
\frac{1}{T} \int_0^T \alpha_0(t)   \diff t = \frac{1}{T} \int_0^T \alpha_1(t)   \diff t = 1.
\]

 We now use Theorem \ref{recast_as_opc} to formulate 
an optimal control problem that allows finding the gain of entrainment. 
Introduce the two-dimensional state $\xi:=\begin{bmatrix} x_{n+1}&x_{n+2}\end{bmatrix}^T$ and let~$u :=\begin{bmatrix} u_0 &u_1 \end{bmatrix}^T$. 
Then the extended system is 
an~$(n+2)$-dimensional nonlinear control system:
\begin{equation}\label{rfm2}
\dot x := \begin{bmatrix} \dot x_1 \\ \dot x_2 \\  \vdots \\ \dot x_n \\  \dot x_{n+1} \\ \dot x_{n+2}   \end{bmatrix}  = f(x)+g(x) u(t) ,
\end{equation}
where
\begin{equation}
f(x)	:= \begin{bmatrix} -\lambda_1 x_1 (1-x_2) \\ \lambda_1 x_1 (1-x_2) - \lambda_2 x_2 (1-x_3) \\ \vdots \\ \lambda_{n-1} x_{n-1} (1-x_{n}) \\ 0 \\ 0    \end{bmatrix} , \quad
g(x):= \begin{bmatrix}   1-x_1  & 0   \\ 0 & 0 \\ \vdots \\ 0 &  -x_{n} \\     \alpha_0(t) & 0    \\  0 & \alpha_1(t)     \end{bmatrix} , \end{equation}
and  the   boundary conditions are: 
\begin{align}\label{boundary}
x_i(T)=x_i(0), i=1, \dots ,n,\;x_{n+1}(0)=0,\;    x_{n+1}(T)=T \bu_0, \;  x_{n+2}(0)=0,\;    x_{n+2}(T)=T \bu_1. 
\end{align}
The optimal control problem is:
\begin{problem} \label{prob:optim}
	Find $u_0, u_1 : [0,T] \to [\ell,L]$ that maximize the cost functional %
	\begin{equation}\label{objective}
	J(u_0,u_1):=\frac{1}{T}\int_0^T \beta(t) u_1(t) x_n(t) \diff t,
	\end{equation}
	subject to the ODE~\eqref{rfm2}, integral constraints~\eqref{eq:oindd},
	and the boundary 
	conditions~\eqref{boundary}, where~$\beta:[0,T] \to \R$ is a given non-negative measurable function. 
\end{problem}

In general, this seems to be a non-trivial problem. Nevertheless, this formulation allows the utilization of both theoretical and
numerical optimal control tools, and we will provide exact results in special cases.

\subsection{Application of the PMP}
As in the previous section, we can write the Hamiltonian \eqref{eq:hamiltonian_general}. In this case,  
\begin{equation} \label{eq:hamiltonian}\mathcal H(u,x,p,p_0)  := p^T(t) \left(f(x(t))+g(x(t)) u(t)  \right) 
+  \tfrac{p_0}{T} \beta(t) u_1(t) x_n(t) ,\end{equation} which   can be written as:
		\begin{align}\label{hamil} 
		\mathcal H &=  \varphi_0(t) u_0(t)+\varphi_1(t)u_1(t),  
		\end{align}
		where  
		\begin{align} \label{e.switching}
		\varphi_0(t)&: = p_1(t) (1-x_1(t))+\alpha_0(t)p_{n+1}(t),\\ \nonumber
		\varphi_1(t)&:= x_n(t) (\tfrac{p_0}{T} \beta(t)-p_n(t)) + \alpha_1(t) p_{n+2}(t) ,
		\end{align}
		are called the \emph{switching functions}.

		\subsubsection{Characterization of  regular arcs} %
		Let $\mathscr X$ be a feasible trajectory. The switching functions $\varphi_0,
		\varphi_1$ play a special role in determining the optimal control. 
	Define the open  set:
		\begin{align*}
		E_r: &= \{ t \in [0,T]\st  \varphi_0(t)\varphi_1(t) \ne 0 \} .
		\end{align*}
 		A \textit{regular arc} is  a restriction  $\mathscr X|_{V}$ for some open subset  $V \subset E_r$.

		Since the Hamiltonian~\eqref{eq:hamiltonian} is linear in the control inputs,   the optimal control is \emph{bang-bang} when the corresponding switching function does not vanish. This is a well-known result in optimal control. We state it  for the sake of completeness.
		\begin{lemma} \label{l.bang} %
			Let $\mathscr X$ be an extremal trajectory. Then for any $t \in E_r$ and~$i\in\{0,1\}$ we have 
			\[ u_i^*(t)=
			\begin{cases}
			    L, \text{ if } \varphi_i(t) > 0 , \\
			\ell,  \text{ if } \varphi_i(t)<0    .
			\end{cases} 
			\]
			This means that at any time $t$ where $\varphi_i(t) \ne 0$, the corresponding $u_i^*(t)$ is a bang-bang control, meaning  that it takes extremal  values.
		\end{lemma}
		\begin{proof}
			We prove the result for~$i=0$. (The proof for~$i=1$ is very similar.)
			Suppose that~$\varphi_0(t)>0$ for some~$t\in [0,T]$. 
			Seeking a contradiction, suppose that $u_0^*(t) <L$. Then,
			\begin{align*}
			\mathcal H(u_0^*(t),u_1^*(t),x^*(t),p^*(t)) &= \varphi_1(t) u_1^*(t)+\varphi_0(t) u_0^*(t)\\&<\varphi_1(t) u_1^*(t)+\varphi_0(t) L \\&= \mathcal H(L,u_1^*(t),x^*(t),p^*(t)), \end{align*}
			and this contradicts~\eqref{eq:pmp_max}. Hence, $u_0^* $ is not optimal. The same argument can be applied when~$\varphi_0(t)<0$.
		\end{proof}
	
	Therefore, unless either of the switching functions vanish on a nonzero measure set, the optimal control   is bang-bang, meaning that it has values in $\{\ell,L\}^2$ for almost all $t$.

	\subsection{The unweighted optimal control problem}
In this section, we consider the unweighted version of Problem~\ref{prob:optim}, that is, the case where~$\alpha_0(t)=\alpha_1(t)=\beta(t)\equiv 1$ for all $t \in [0,T]$.  

\subsubsection{Constant controls  satisfy  the PMP}
We first show that the constant controls
satisfy  the necessary conditions for optimality.  
\begin{theorem} \label{t.constant}
Consider the unweighted Problem \ref{prob:optim}. The constant controls~$u_0(t)\equiv \bu_0$, $u_1(t) \equiv \bu_1$ satisfy Proposition \ref{p.pmp} (the PMP) with the corresponding switching functions  identically zero. 
\end{theorem}
\begin{proof}
	Let~$z:=\begin{bmatrix}p_1&\dots&p_n \end{bmatrix}^T$, i.e. the first~$n$ entries of the adjoint state. 
	Eq.~\eqref{pdot} yields
	\be\label{eq:pnaur}
	\dot z= -J^T(  x,u ) z -b u_1,
	\ee
	where $J$ is the Jacobian of the RFM \eqref{rfm1} w.r.t.~$x$, 
	and~$b:=\begin{bmatrix}0 & \dots &0 &  p_0    /T \end{bmatrix}^T$.
	Also,~$\dot p_{n+1}(t)=\dot p_{n+2}(t)\equiv 0$.

	It has been shown in~\cite{margaliot2012stability}  that the RFM with constant rates  admits a unique GAS steady state in~$(0,1)^n$. Hence, 
 every solution  of~\eqref{rfm1} with $u_0(t)\equiv\bu_0, u_1(t)\equiv \bu_1$,
converges to a point~$\bar x=\begin{bmatrix}\bar x_1& \bar x_2&\dots&\bar x_n\end{bmatrix}^T \in (0,1)^n$.

	  It was shown
	in~\cite{entrainment} that if~$M$ is any compact subset of~$(0,1)^n$ then
	there exists a matrix measure~$\mu:\R^{n\times n} \to \R$ such that~$\mu(J(x,u))<0$ for all~$x\in M, u\ge 0$.
	This implies in particular that all the eigenvalues of~$J(x,u)$ have a 
	negative real part~\cite{sontag_cotraction_tutorial}, so~$J(x,u)$ is nonsingular for each $x \in (0,1)^n , u\ge 0$. Hence, so is~$J^T(\bar x,\bar u)$. Let~$\bar z :=- (J^T(\bx,\bu))^{-1} b \bu_1$.  	Let $\bu:=\begin{bmatrix} \bu_0&\bu_1\end{bmatrix}^T$. We now show that for
	\[
	u(t)\equiv \bu,\; x(t)\equiv \bar x, \; p_0=T, \; p(t) \equiv \begin{bmatrix} \bar z &
  -\bar p_1 (1-\bar x_1)/\bu_0  &  -\bar x_n (1-\bar p_n)/\bu_1 \end{bmatrix}^T ,	
	\]
	all the conditions in the PMP hold. First note that the boundary conditions~\eqref{boundary} all hold. 
	Eq.~\eqref{eq:pnaur} holds by the definition of~$\bar p$. 
	The switching functions~\eqref{e.switching} satisfy 
	 $\varphi_0(t)\equiv\varphi_1(t)\equiv 0$. 
	Eq.~\eqref{hamil} implies that~$\mathcal H\equiv0$ and that~\eqref{eq:pmp_max} trivially  holds.  
\end{proof}

We have shown that   constant controls satisfy the necessary conditions of the PMP.  In other words,   constant controls are always extremal solutions.

In the next subsection, we show that for~$n=1$
  constant controls are the \emph{only} controls that satisfy 
  the~PMP.
\subsubsection{Extremal analysis of the  one-dimensional unweighted problem}
  In this subsection, we study the following system:
\begin{align}  \nonumber
\dot x_1(t) &= u_0(t) (1-x_1(t)) - u_1(t) x_1(t), \\\label{rfm1d}
\dot x_2(t) &= u_0(t), \\ \nonumber
\dot x_3(t) &= u_1(t).
\end{align} 
The controls~$u_0$ [$u_1$] 
represent     time-varying initiation [exit] rates in an  RFM with~$n=1$. 
	Even though the PMP provides a general approach for addressing optimal control problems, it seldom leads to  a  full characterization of extremal solutions, especially in the case of   multiple inputs. 
We will show that this is possible for the system \eqref{rfm1d}: a detailed analysis using the  PMP shows that any
extremal trajectory corresponds to a  constant~$x_1(t)$. Since each control input takes values in a compact and convex set, the optimal control problem always has a solution. Thus, there is no gain of entrainment. 
 This shows that the PMP is a viable approach for handling such problems and lays the ground for future generalization to higher dimensional cases.

	 \begin{theorem}\label{th.extremal}
	 	Let $\mathscr X$ be an extremal trajectory for Problem \ref{prob:optim} with the system \eqref{rfm1d}. 
		Then
		\be\label{eq:xgtre}
		x_1^*(t)\equiv \frac{\bar {u}_0 }{\bu_0+ \bu_1}=\frac 1{1+(\bu_1/\bu_0)} \text{ for all } t \in [0,T].
		\ee
	 \end{theorem}
\noindent The proof is given in the Appendix.

The PMP   immediately yields the following.
\begin{theorem}\label{theorem2in}
	Fix  $0<\ell<L$, and let the admissible controls ~$u_0,u_1$ take values in~$ [\ell,L]$,  with
	given averages~$\bu_0 ,\bu_1 \in(\ell,L) $,
	the optimal objective for Problem~\ref{prob:optim} and 
	system \eqref{rfm1d} is
	\begin{equation}\label{optimal_cost}
	J^*= \frac{\bu_0 }{\bu_0 + \bu_1 }.
	\end{equation}
	The optimal   trajectory is~$x_1^*(t)\equiv \bu_0 / (\bu_0 + \bu_1)$, a constant, 
	and it can be achieved by the constant  inputs:
	\begin{equation}\label{eq:optimal_u} 
	u_0^*(t) \equiv  \bu_0, u_1^*(t)\equiv \bu_1.
	\end{equation}
\end{theorem}
\begin{remark}
	The control inputs~$u_0(t), u_1(t)$ that achieve  the optimal cost are not unique.  
	Indeed,  it is clear that the    optimal solution in~\eqref{eq:xgtre}
 depends only on the ratio~$\bu_1/\bu_0$. For instance, $u_0^{**}(t) \equiv  \bu_0 \rho(t), u_1^{**}(t)\equiv \bu_1 \rho(t)$ is
 also an optimal solution for any function~$\rho$ such that
 $ \bu_0 \rho(t)  \subset [\ell,L]$ and 
 $ \bu_1 \rho(t)  \subset [\ell,L]$
 for all $t\in[0,T]$.
\end{remark}

\begin{remark}
Theorem \ref{theorem2in} can be also proven  via a more direct approach, motivated by an idea from~\cite{katriel20}. This alternative proof is given in the Appendix. 
\end{remark}

\subsubsection{The   unweighted problem for the RFM with~$n=2$ and a single input}
We now study Problem~\ref{prob:optim} for 
 an RFM with~$n=2$ and a single control~$u_0(t)$ 
as the initiation  rate, i.e. 
\begin{align} \nonumber
\dot x_1 &=  u_0 (1-x_1)  - \lambda_1 x_1 (1-x_2 ) , \\ \label{2drfm}
\dot x_2 &= \lambda_1 (1-x_2) x_1 -  \lambda_2 x_2 ,\\ \nonumber
\dot x_3 &= u_0(t),
\end{align}
subject to the boundary conditions~\eqref{boundary}. Here our goal it to maximize the average value  of the output rate~$\lambda_2 x_2 $.

More formally, we aim at solving the following problem:
\begin{problem} \label{prob:optim_1in}
Let $\ell,\bu_0,L$ be given such that $0<\ell<\bu_0<L$. Find $u_0 : [0,T] \to [\ell,L]$ that maximizes the cost functional %
$
	J(u_0):=\frac{1}{T}\int_0^T  \lambda_2 x_2(t) \diff t,
$
that is, the average production rate 
subject to the ODE~\eqref{2drfm}, integral constraint~$\frac 1T \int_0^T u_0(t)=\bu_0$,
	and the boundary 
	conditions $x_1(0)=x_1(T), x_2(0)=x_2(T)$. 
\end{problem}

The next result 
shows   that here as well a constant control is optimal, that is, there is no gain of entrainment.  
  
\begin{theorem}\label{result2d}
	 Fix~$0<\ell<L$ and~$\bu_0 \in (\ell,L)$,  
	the objective function for Problem~\ref{prob:optim_1in}  with the system \eqref{2drfm} is maximized by the constant control~$u_0^*(t)\equiv \bar u_0$.
\end{theorem}
\begin{proof}
The first equation in~\eqref{2drfm} is the first equation~\eqref{rfm1d} for~$u_1(t)=\lambda_1(1-x_2)$, so~$\bu_1=\lambda_1(1-\bar x_2)$, and the output rate is~$\lambda_1  (1-x_2) x_1$. 
(Note that here  $x_2$ is the second state-variable in~\eqref{2drfm},  and  not  the integral of $u_0$ as in~\eqref{rfm1d}). 
	By Theorem~\ref{theorem2in}, we have $\lambda_1 \overline{ u_1 x_1} \le (\lambda_1\bu_1 \bu_0)/(\bu_0+\lambda_1\bu_1 )$. Hence,
\begin{equation}\label{th7_p}	\lambda_1 \overline{ u_1 x_1} =\lambda_1 \overline{(1-x_2) x_1 } \le \frac{\lambda_1\bu_1\bu_0}{\bu_0+\lambda_1\bu_1 } =\frac{ \lambda_1   (1-\bx_2) \bu_0}{ \bu_0+\lambda_1(1-\bx_2)}.\end{equation}
Integrating \eqref{2drfm} we get \[0=x_2(T)-x_2(0)=\int_0^T \dot x_2(t) \diff t = \lambda_1 \int_0^T (1-x_2(t))x_1(t) \diff t - \int_0^T \lambda_2 x_2(t) \diff t.\] Hence, we get that $   \lambda_1 \overline{(1-x_2) x_1 }=\lambda_2 \bx_2    $. Substituting in \eqref{th7_p}, we get
\be\label{eq:quntmax}
\lambda_2  \bx_2 \le  \frac{ \lambda_1 \bu_0 (1-\bx_2)}{\lambda_1(1-\bx_2)+\bu_0}.
\ee
The left-hand side here is the quantity that we are trying to maximize.  
Rearranging gives: 
\be\label{eq:fxwe}
f(\bar x_2)\ge 0,
\ee
 where $f(s):= s^2 - \left( 1 + \bu_0(\tfrac 1{\lambda_1}+ \tfrac 1{\lambda_2}) \right ) s + \tfrac{\bu_0}{\lambda_2}$.  
Let~$p,q$ denote the roots of~$f(s)$. Then~$f(s)=(s-p)(s-q)$ gives
\begin{align}\label{eq:posaee}
                        1+\bu_0 (\tfrac 1{\lambda_1}+ \tfrac 1{\lambda_2}) =p+q , \nonumber \\
                          \tfrac{ \bu_0} {\lambda_2}  =p q.
\end{align}

Recall that the RFM~\eqref{2drfm} with the constant control~$u_0(t)\equiv \bar u_0$  admits a \emph{unique} steady 
state~$\bar e\in(0,1)^2$~\cite{margaliot2012stability}. It is straightforward 
to show that~$f(\bar e_2)=0$. We may assume that~$p=\bar e_2$, so~$p\in(0,1)$.
Then~\eqref{eq:posaee} implies that~$q$ is real and~$q>1$.  
The quadratic inequality \eqref{eq:fxwe} implies that
either~$\bar x_2\le p<1$ or $\bar x_2 \ge q > 1$. Since~$x_2(t) \in [0,1]$ for all~$t$, then the second inequality can be ignored and we
have that the maximal (feasible) value of~$\bar x_2$ is~$\bar x_2^*=p=e_2$. 
Obviously, this is attained  for the constant control~$ u_0(t)\equiv \bar u_0$.  
\end{proof}

 		\subsection{Gain of entrainment with   time-varying weight  functions}
Analyzing the case  with time-varying weights is challenging, but it is highly relevant to applications since resources may be
 allocated differently  during the period. 
In this subsection, we show that, even for the above examples,
once the weighting functions become time-varying,
constant inputs may no longer be optimal.

 We consider the special case of~\eqref{rfm2} with $n=1$ and a single input $u_0(t)$ 
as the initiation rate, i.e. 
\begin{align} \label{tv} \dot x_1(t) &= u_0 (t)(1-x_1(t)) -  \lambda_1 x_1(t), \\ 
\dot x_2(t)&=u_0(t).\nonumber 
 \end{align}
Our goal now is to maximize
	$J(u_0):=\frac{\lambda_1}{T}\int_0^T \beta(t)    x_1(t) \diff t$,
	subject to the boundary conditions~\eqref{boundary}, and where
the  weight  function~$\beta$ is differentiable and  satisfies  $\beta(t)>0$ for all~$t \in [0,T]$. 
Without loss of generality, we assume that $T=1$.

\begin{proposition} \label{pr.tv}
Suppose that~$u_0^*$ is an optimal control. 
Then for almost all~$t\in[0,1]$ we have that  either~$u_0^*(t)\in \{\ell, L\}$ or:
\begin{equation}\label{e.tv_singular} u_0^*(t) = c \sqrt{\beta(t)} - \lambda_1 + \frac{\dot \beta(t)}{2\beta(t)}, \end{equation}
for some constant $c$. Furthermore, if $u_0^*$ satisfies~\eqref{e.tv_singular} 
 for all $t \in [0,1]$, then \begin{equation}\label{constant}  c = \frac{ \bu_0+\lambda_1 - \tfrac 12 \log \tfrac{\beta(1)}{\beta(0)}}{\int_0^1 \sqrt{\beta(t)} \diff t }. \end{equation}
\end{proposition}

Note that if~$\bar u_0 \in (\ell,L)$  
then this  implies that the constant control~$u_0(t)\equiv \bar u_0$ cannot be optimal. 
If~\eqref{e.tv_singular}  does not hold for any~$t$ (e.g. when the right-hand side of~\eqref{e.tv_singular} takes values that are not in~$[\ell,L]$) then the optimal control is bang-bang. 
 
As a specific example 
take
\[
\ell=0.001,\;L=10, \bar u_0=2, \lambda_1=1,
\]
and the  weight  function
\[
 \beta(t)= e^{ -\rho \left (t-  (T/2) \right )^2},
\]
 with~$\rho=100$ and~$T=1$. 
In the context of the RFM, this would  represent  the case  where
 it is required to highly  expresses a specific protein  near 
 the middle of every cycle, rather than having a uniform level of production along the entire  period.

The constant input $u_0(t)\equiv   2$  yields a steady-state trajectory~$x(t)\equiv   2  /3 $. The corresponding value of the objective function 
 is
\be\label{eq:juop}
J(u_0)=\frac{2}{3}\int_0^1 \beta(t)     \diff t =0.118164.
\ee

Our optimal control formulation   can be used to solve the problem \emph{numerically} using
optimal control packages like~\cite{gpops}. 
The result is a 
 three-arc bang-bang control  
\[
u^*(t)=\begin{cases}  
\ell,& t\in[0,t_1) ,\\ 
L,& t\in[t_1,t_1+\Delta) ,\\ 
\ell,& t\in[t_1+\Delta,1],
\end{cases} 
\]
where~$\Delta:=(\bu_0-\ell)/(L-\ell)=0.19992$, and $t_1=0.27335$.
 The corresponding periodic solution satisfies~$x(1)=x(0)= 0.50251$.
This  achieves a cost~$J(u^*)=0.141183$,
 which is roughly 20\% better than~\eqref{eq:juop}. 

\begin{figure}
 \centering
\includegraphics[width=0.65\textwidth]{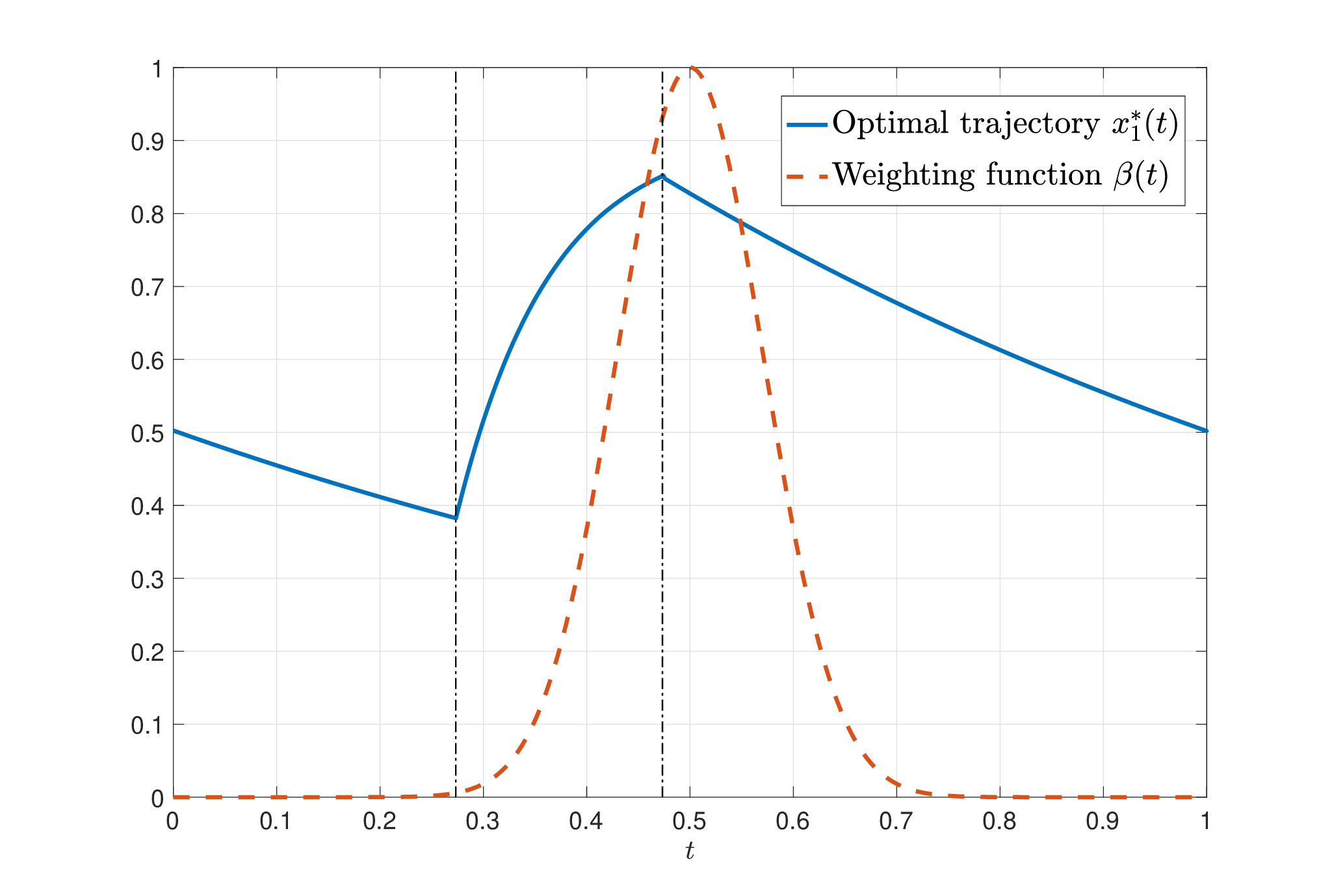} 
\caption{The optimal trajectory for maximizing a weighted throughput is \emph{non-constant}. The plot shows the weighting function $\beta(t)$ with the corresponding  optimal periodic solution~$x_1^*(t)$. The vertical dashed lines denote  the switching points of the three-arc bang-bang controller. Note that the high bang control is switched on when the value of the weight  function becomes non-negligible.}\label{f.TV}
 \end{figure}

Fig.~\ref{f.TV} depicts the approximate optimal bang-bang control   (computed using \cite{gpops} and a bisection procedure)
and the  resulting periodic solution~$x_1^*(t)$. The weight function~$\beta(t)$ is also shown. 
The maximal value of~$\beta$ is achieved at~$t=1/2$. The 
  optimal control switches to  the maximal  value~$L=10$
  before the peak time of~$\beta$. 
	This makes sense  as it guarantees that~$x_1^*(t)$ 
	will have large values when the 
	weighting of~$x_1^*(t)$ in the cost function is large.

\section{Gain of entrainment in generalized occupancy models}\label{sec:entrain_GOM}
In this section, we consider  general cascades like those in in Figure~\ref{f.cascade}.  Recall that under the conditions~\eqref{eq:adftp},  the $n$-dimensional RFM can be approximated by \eqref{eq:rfm_l}, which  has the form  in Figure~\ref{f.cascade}-(a). We consider here this approximated system.

We  first state the following  result which is well-known in the theory of linear time-invariant systems (see also \cite{RFM_IEEE_CL}). We include the proof for completeness.%
\begin{proposition} \label{prop:9}
	Consider a    singe-input-single-output   linear system: $\dot z=A z+b w $, $y=c^T z$, where $z \in \mathbb R^n$, and~$A$ is Hurwitz.
	Let~$w$ be a bounded measurable  input which is~$T$-periodic.
	Then~$y$ converges to a steady-state~$T$-periodic solution~$y_w$, and 
	\[
	\int_0^T y_w(t) \diff t
	= H(0) \int_0^T  w(t) \diff t  ,
	\]
	where $H(s):=c^T(sI-A)^{-1}b$ is the  transfer function of the linear system.
\end{proposition}

\begin{proof}
	Since $w$ is measurable and bounded,~$w \in L_2 ([0,T])$. Hence, it can be written as a Fourier series~$
	w(t)=\bar w+ \sum_i a_i \sin (\omega_i  t +\phi_i )$.
	The output of the linear system
	converges to the steady-state periodic solution
	\begin{align*}
	y_w(t)&=    H(0) \bar w  + 
	\sum_i  |H(j\omega_i)| a_i  \sin (\omega_i  t +\phi_i  +\angle (H(j\omega_i) )).
	\end{align*}
 	Each sinusoid in the expansion has period~$T$, so
	$ \int_0^T y_w(t) \diff t = T  H(0) \bar w .$
\end{proof}

Combining  Proposition~\ref{prop:9}  with our results on the gain of entrainment in certain bottleneck models  
yields the following result. 

\begin{theorem}\label{theorem_LS}
	Consider the nonlinear system depicted in either Figure \ref{f.cascade}-(a) or Figure \ref{f.cascade}-(b) with $A$  Metzler, and~$b,c \in \mathbb R_{+}^n$. Let $u_0(t), u_1(t)$ be $T$-periodic non-negative control signals.	
	For any~$0<\ell<L$ and any~$\bu_0,\bu_1  \in (\ell,L)$,  
	consider	the functional	\[
	J(u_0,u_1):=\frac{1}{T}\int_0^T w_2(t) \diff t,
	\]
	where $w_2(t)$ is the steady-state $T$-periodic output signal. Then  the constant  controls:
	\begin{equation}\label{optimal_u} 
	u_0^*(t) \equiv  \bu_0, \  u_1^*(t)\equiv \bu_1.
	\end{equation}
	maximize $J$.
\end{theorem}
\begin{proof}
Consider the system depicted in Figure \ref{f.cascade}-(a).  
 Fix admissible~$T$-periodic controls~$u_0(t)$, $u_1(t)$. 
Denote the corresponding steady-state  average values of~$w_1(t)$ and~$y(t)$ 
by~$\bar w_1$ and~$\bar y$.  Obviously,~$\bar w_2=\bar y$.
By Proposition \ref{prop:9}, $\bar w_2=H(0) \bar w_1$, where~$H$ is the transfer function of the linear system. Since~$A$ is Metzler
and~$b,c \in \mathbb R_{+}^n$, the trajectories of the linear system are positive. Thus, maximizing~$J$ is equivalent to maximizing~$\bar w_1$.  Theorem \ref{theorem2in}  implies that  the constant controls~$u_0^*(t)\equiv \bu_0$ and~$u_1^*(t)\equiv \bu_1$ maximize~$\bar w_1$. 
	The system in   Figure \ref{f.cascade}-(b) can be treated similarly.
\end{proof}

\begin{figure}
	\centering
	\includegraphics[width=0.75\textwidth]{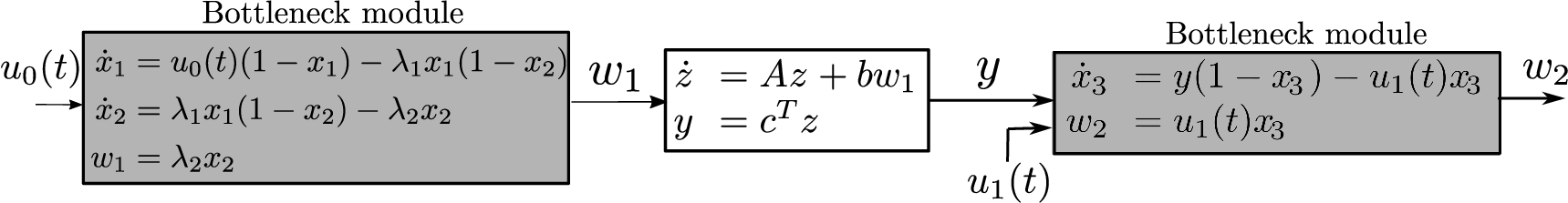}
	\caption{A generalized occupancy model with a 2D bottleneck entrance. The controls are $u_0(t),u_1(t)$ which are scalar functions. We have $x_1, x_2, x_3, w_1, w_2, y \in \mathbb R_{+}$, $z \in \mathbb R^n, A \in \mathbb R^{n \times n}, b,c \in \mathbb R_{+}^n$. The linear system block is assumed to be positive and Hurwtiz.}
	\label{f.2d_sandwich}
\end{figure}
\begin{remark}
The same result holds if 
	the single-input modules in Figure~\ref{f.cascade} are
	replaced by the single-input  RFM with $n=2$ in~\eqref{2drfm} as shown in Figure \ref{f.2d_sandwich}.
\end{remark}

 \section{Discussion}\label{sec:conc}
Entrainment is an interesting and
important property of dynamical systems.  
It allows systems to develop an ``internal clock'' that synchronizes to periodic excitations like the 24h solar day.
Such clocks are important in biology, as they allow organisms to adequately respond to periodic processes like the   solar day and the cell cycle division process. They are also essential for synthetic biology, as a  common clock is an important ingredient in building synthetic biology circuits that include several modules working in synchrony. 

Here, 
we considered an additional 
  qualitative property  called the \emph{gain of entrainment}.
	This   measures the advantage, if any,  of using 
a periodic control vs. an ``equivalent'' constant control for maximizing the  average throughput. We showed how this problem can be cast as an optimal control problem. 
This allows using the sophisticated analytical and numerical tools developed for solving optimal control problems to determine the gain of entrainment. 

 We have shown that, perhaps surprisingly,
 there is no gain of entrainment in  a class of systems relevant to biology and traffic applications. In other words, in these systems non-constant periodic controls are no better than constant controls.  The optimality of constant controls fails to hold if we allow  \emph{time-varying} objective functionals. Hence, this suggests  the possibility that the observation of non-constant periodic signals in biological systems is correlated with the maximization of the throughput with varying weights along ecah period or cycle.

An interesting research direction is to generalize these results to models such as the nonlinear $n$ site RFM  with $(n+1)$ time-periodic control inputs.  
%
%

%
%
 \begin{comment}

hold off;

fun = @(t,v) exp(  t+sin(v*t)/2/v ) ;
xt= @(t,v,d) (exp( -t - sin(v*t)/2/v ).*( d + exp(  t+sin(v*t)/2/v ))) ;

for w=0.5:.01:5
    
     T=2*pi/w;     
     c=exp(-T)* integral(@(t)fun(t,w),0,T)/( 1-exp(-T) );
          
	 val=(1/T)*integral(@(t) xt(t,w,c),0,T)-1;
     
	 plot(w,val,'.k'); hold on;
end
h=xlabel( '$\omega$','interpreter', 'latex' );
 %
  set(h,'FontSize',18);
grid on;
%

\end{comment}
 
%
	
	\section*{Appendix: Additional Proofs}
	\subsection*{A.1 \quad Proof of Proposition \ref{p.pmp}}
		Most of the statements here are   the
		standard~PMP. 
		We only need to prove  the transversality condition~\eqref{eq:trans}.%
		
		Pick a set~$S\subseteq  \R^{2(n+m)}$, 
		and suppose that the state must satisfy the  constraint~$ \begin{bmatrix} x(0)^T & x(T)^T \end{bmatrix}^T \in S$. 
		Then 
		the corresponding transversality condition~\cite{LeeMarkus}
		is
		\[ 
		\begin{bmatrix} p(0) \\ -p(T) \end{bmatrix} \bot \cT_{  \text{\tiny $\begin{bmatrix} x(0) \\ x(T) \end{bmatrix}$} } S, %
		\]
		where $\cT_{z } S$ denotes the  tangent space 
		of~$S$ at~$z $.
		In our case, \eqref{boundary_general} gives
		\[
		S = \{ z \in \R^{2n+2m}\st  [z_1,..,z_n]^T = [z_{n+m+1},..,z_{2n+m}]^T, [z_{n+1},..,z_{n+2}]^T=0, 
		[z_{2n+m+1},...,z_{2n+2m}]^T=T q\}.
		\]
		Hence, $T_z S=\mbox{span}\{ v^1,...,v^n\}$, 
		where  $v^i$ is the vector with one at entries~$i$ and~$(i+n+m)$, and zero elsewhere. 
		Therefore, it is necessary that $p_i^*(0)=p_i^*(T),~ i=1,\dots ,n$.	 \hfill $\blacksquare$

	\subsection*{A.2 \quad Proof of Theorem \ref{th.extremal}}
	The proof, based on the analysis of extremals,  is divided into a sequence of Lemmas. 	 For a set~$A \subset \mathbb R$,
	$\mu(A)$ denotes its Lebesgue measure. The set of accumulation 
	points of~$A$ is denoted by~$A'$. 
	For $x \in \mathbb R$,  $\{x\}+A:=\{ x+a \st  a \in A \}$. 
	
	Recall that we consider \eqref{rfm1d} , so
	\[\mathcal H = p_1 ( u_0 (1-x_1)-u_1 x_1) +p_2 u_0 + p_3 u_1 +  \tfrac{p_0}T u_1 x_1. \]
	
	\begin{lemma}\label{l.pode}
The   adjoint variables $p_i^*$ satisfy:
\begin{align}\label{p_ode}
\dot p_1^*(t) &= - (u_0(t)+u_1(t)) p_1^*(t) - u_1 (t), \\ \label{e.p2p3} 
p_2^*(t)&\equiv p_2^*(0),  \nonumber \\
p_3^*(t)&\equiv p_3^*(0),\nonumber 
\end{align}
with the boundary condition $p_1^*(0)=p_1^*(T)$.
	\end{lemma}
\begin{proof}
	We first      show that we can take~$p_0^* =T$.
	Assume that~$p_0^*=0$. Then~\eqref{pdot}
	yields
	$	\dot  p_1^* = (u_0^*(t)+u_1^*(t) ) p_1^*.$ Integrating over $[0,T]$, we get 
	\[
	\log |p_1^*(T)|-\log |p_1^*(0)|=\int_0^T (u_0^*(t)+u_1^*(t) ) \diff t = T(\bu_0+\bu_1) . 
	\]
	By the 
	transversality condition, we know that~$ p_1^*(0)= p_1^*(T)$, which implies that $\bu_0+\bu_1=0$. This  is a contradiction, so we conclude that~$p_0^*> 0$, and by scaling the objective function we may take~$p_0^*=T $. Now~\eqref{p_ode} follows from calculating the partial derivatives in~\eqref{pdot}.  
	\end{proof}

	\subsubsection*{Analysis of the switching functions}
	Using  Lemma \ref{l.pode}, the switching functions in our case are:
	\begin{align} \label{phi0}
	\varphi_0(t)& = p_1(t) (1-x_1(t))+p_{2}(0), \\ \label{phi1}
	\varphi_1(t)&= x_1(t) (1-p_1(t)) +  p_{3}(0) . 
	\end{align}
	
	Given an extremal trajectory $\mathscr X$, let
	\begin{align*}
	E_+^i:=\{ t \in [0,T] \st  \varphi_i (t)>0\} , \\
	E_-^i :=\{ t \in [0,T] \st \varphi_i (t)<0\},\\
	E_0^i : =\{ t \in [0,T] \st   \varphi_i (t)=0\},
	\end{align*}
	where $i=0,1$. Note that $E_+^i, E_-^i$, $i=0,1$, are open relative to $[0,T]$, and $E_0^i, i=0,1$, are closed. In particular, all these  sets are Lebesgue measurable. 
	
A calculation gives
	\begin{align}
	\dot\varphi_0(t) &= u_1(t) ( p_1(t) - (1-x_1(t))), \label{dotphi0}\\
	\dot\varphi_1(t) &= u_0(t) (1-x_1(t)-p_1(t)). \label{dotphi1}
	\end{align}
	\begin{remark} \label{r.phi_der} The functions $\varphi_0,\varphi_1$ are absolutely continuous. Hence, they are differentiable almost everywhere and have bounded derivatives. This implies that both $\varphi_0,\varphi_1$ are Lipschitz continuous. Also,  since the controls are positive, \eqref{dotphi0},\eqref{dotphi1} imply that $\sgn(\dot\varphi_0(t))=-\sgn(\dot\varphi_1(t))$ whenever  $\varphi_0,\varphi_1$ are both differentiable.\end{remark}
	\subsubsection*{Characterization of singular arcs}

	In this subsection, we are interested in  the case 
	where~$\mu(E_0^i)>0$ for either $i=0$ or $i=1$. 
	Let  \[E_s:=\{t \in [0,T] \st \varphi_0(t)\varphi_1(t)= 0 \} = E_0^0 \cup E_0^1.\]
	Let $\mathscr X$ be an extremal trajectory. We call any restriction  of $\mathscr X$ to any nonzero-measure subset of $E_s$  a \emph{singular arc}.
	
	The following Lemmas characterize the   behavior on singular arcs.

	\begin{lemma}\label{l.s1}
		Let $\mathscr X$ be an extremal trajectory, and assume that~$\mu(E_0^i)>0$ for some $i\in\{0,1\}$.
		Then there exists %
		$c_i\in (0,1)$ such that 
		\[
		x_1^*(t)=  c_i \text{ for almost all} \ t \in E_0^i.
		\] 
		Furthermore, $\dot x_1^*(t)=0$ for almost all $t \in E_0^i$ and the two inputs   satisfy:
		\begin{equation}\label{inputs_coupled}\left ( \frac 1{c_i}-1 \right )	u_0^*(t) =  u_1^*(t), \;\mbox{for almost all} \ t \in E_0^i.
		\end{equation}
	\end{lemma}
	\begin{proof}
		Let $ E_0^{i'} \subseteq E_0^i$ denote  the set of accumulation points of~$E_0^{i}$.
		Note that~$\mu(E_0^{i'})=\mu(E_0^i)$, since $E_0^{i}\setminus E_0^{i'}$ is the set of 
		isolated points of~$E_0^i$ which is countable, and hence has measure zero.  Let $F_i:=\{t \in [0,T] \st  \dot \varphi_i^*(t) \ \mbox{exists} \} \cap E_0^{i'}$. Then~$\mu(F_i)=\mu(E_0^i)$  since $\varphi_i^*$ is differentiable a.e.
		
		Fix~$t\in F_i$. By definition we have 
		$\varphi_i^*(t)=0$. We show that $\dot\varphi_i^*(t)=0$ as well. Since~$t$ is an accumulation point, $\exists \{t_k\} _{k=1}^\infty \subset E_0^i$ such that $t_k \to t$. Since $\varphi_i^*$ is differentiable at~$t$,   $\dot \varphi_i^*(t) = \lim_{k \to \infty} (\varphi_i^*(t_k) -\varphi_i^*(t))/(t_k-t)=0$. We consider the case~$i=0$ 
		(the proof when~$i=1$ is very similar). 
			Using \eqref{phi0} and~\eqref{dotphi0}	the equations~$\varphi_0^*(t)=\dot\varphi_0^*(t)=0$ yield~$p_1^*(t) (1-x_1^*(t)) = -p_2^*(0)$, and $p_1^*(t) = 1-x_1^*(t)$. Since~$x^*_1(t) \in (0,1)$,
		   $   p_2^*(0)<0$ and 
				\begin{equation}\label{eq:xqpo}
		x_1^*(t)= 1-p_1^*(t) = c_0, \; \mbox{for all} \; t \in F_0.
		\end{equation}
		where~$c_0 := 1-\sqrt{-p_2^*(0)}$.
		
				Let $F_0^{'}$ be the set of accumulation points of $F_0$.
				Then~$\mu(F_0^{'})=\mu(E_0^0)$. Fix~$t \in F_0^{'}$. Hence  $\exists \{t_k\} _{k=1}^\infty \subset F_i$ such that $t_k \to t$. Since $x_1^*(t_k)=c_0$ for all $k$,    $\dot x_1^*(t)=0$. 
				Substituting this in~\eqref{rfm1d} proves~\eqref{inputs_coupled}. 
	\end{proof}

	The next result shows  that if an extremal trajectory $\mathscr X$ has $x_1^*$ identically constant, then it  satisfies~\eqref{eq:xgtre},
	  and~$\mathscr X$   consists entirely of singular arcs.
	\begin{lemma} \label{l.constant}
		Let $\mathscr X$ be an extremal trajectory. If $x_1^*(t)$ is identically constant on~$[0,T]$, then 
		\be\label{eq:x1trp}
		x_1^*(t)\equiv \bu_0/(\bu_1+\bu_0).
		\ee
		Furthermore, $\varphi_0^*(t)\varphi_1^*(t)\equiv 0$, i.e. $E_s=[0,T]$.
	\end{lemma}
	\begin{proof}
		By assumption,  there exists~$c\in(0,1)$ such that~$c\equiv x_1^*(t)$.  Substituting this in \eqref{rfm2} yields~$\dot x_1^* = u_0^*(t)-(u_1^*(t)+u_0^*(t)) c \equiv 0$, and integrating over~$[0,T]$ proves~\eqref{eq:x1trp}.  We also find that~$ u_1^*(t)\equiv (\tfrac 1c -1 ) u_0^*(t) $, and substituting this in~\eqref{p_ode} gives~$\dot p_1^*(t)=u_1^*(t) (-1-\tfrac1{1-c}p_1^*(t))$. Combining this with  the boundary condition~\eqref{eq:trans}
		gives~$p_1^*(t)\equiv c-1$. 
		 Eqs.~\eqref{phi0} and~\eqref{phi1}  give
		\begin{align}\label{eq:pkytre}
		 \phi_0(t)&=  (c-1) (1-c)+p_2^*(0)  ,\nonumber \\
		 \phi_1(t)&=   c ( -c)+p_3^*(0)   ,
		\end{align}
		for all~$t\in[0,T]$. 

		To prove that $E_r = \emptyset$, we first assume that
		that~$\mu(E_s)=0$. Then~$E_r=[0,T]$ (a.e.), and
	Eq.~\eqref{eq:pkytre} implies that both switching functions are constant and not zero, by the definition of $E_r$.  Lemma~\ref{l.bang} implies that every~$u_i^*$ 
	is constant, i.e., either $u_i^*(t)\equiv\ell$ or $u_i^*(t)\equiv L$.
 This contradicts the fact that~$\bar u_i^*=\frac{1}{T}\int_0^T u_i^*(s)\diff s  \in (\ell,L)$.
Thus,~$E_s=[0,T]$. We conclude that~$\mu(E_s)>0$. 
		Pick~$\tau\in E_s$. Then~$\phi_0(\tau)=\phi_1(\tau)=0$ and~\eqref{eq:pkytre}
		implies that~$\phi_0(t)=\phi_1(t)=0$ for all~$t$, so~$E_s=[0,T]$.
	\end{proof}
	
	\subsubsection*{Inadmissibility of Regular Arcs}
	
	In the previous subsection we  decomposed an extremal trajectory into regular and singular arcs. On the regular arcs, the control is bang-bang.  On the singular arcs, the controls  satisfy~\eqref{inputs_coupled} and the state must be constant almost everywhere, so~$\dot x_1^*(t)=0$ a.e. %
	In general, an extremal trajectory can consist of an arbitrary patching of regular and singular arcs. In this section, we consider  the admissibility of regular arcs. 

\begin{figure}
		\centering
		\includegraphics[width=0.8\textwidth]{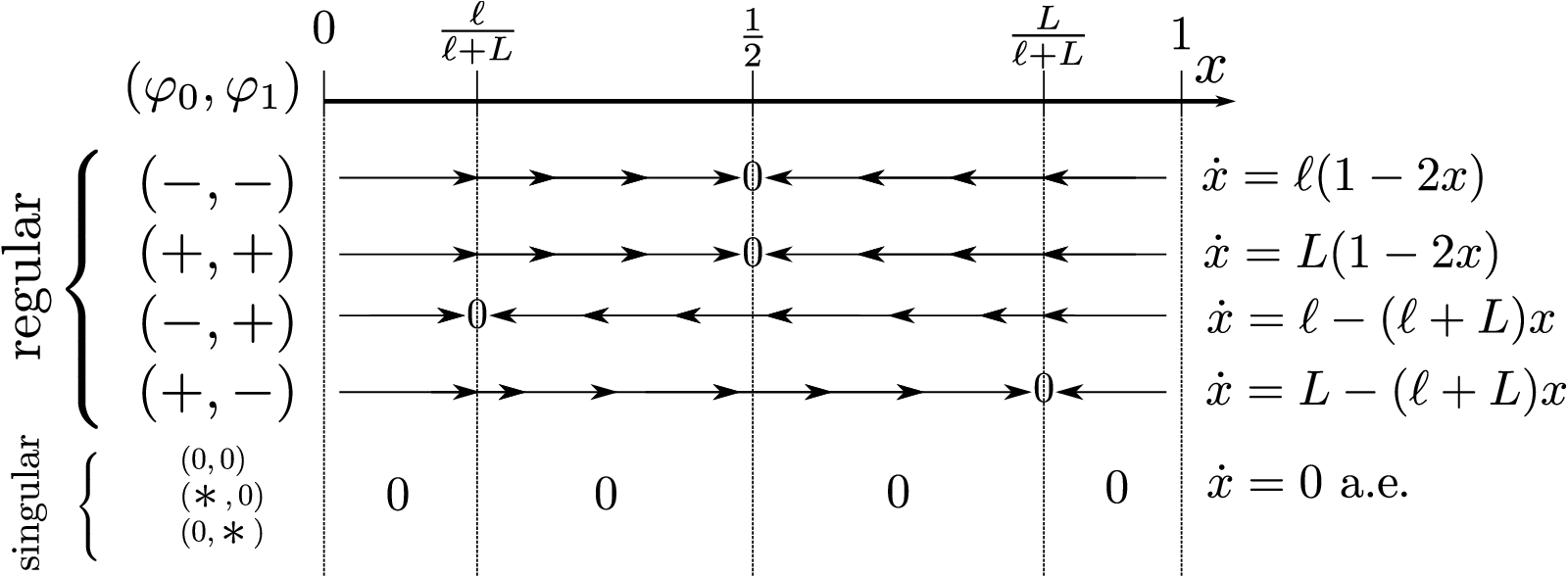}
		\caption{ The equation for~$\dot x $ 
		and the directions of the dynamics  as a function of~$x\in(0,1)$
		for all possible arcs in an extremal trajectory (see Lemmas~\ref{l.bang} and~\ref{l.s1}).
		A circle on an axis describes an equilibrium point of the dynamics. 
		The same diagram   holds for~$p $, but with all the arrow directions reversed. }
		\label{f.signs}
	\end{figure}

	To simplify presentation, we   write $x$ and~$p$ instead of~$x_1$ and~$p_1$ from here on. 
	Furthermore, we   denote extremal trajectories by~$x^*,p^*$.
	Figure~\ref{f.signs} depicts the dynamics of~$x$ based on Lemmas~\ref{l.bang} and~\ref{l.s1}.
	
	\begin{lemma} \label{l.bounds}
		Let $\mathscr X$ be an extremal trajectory. Then $x(t), p(t)  \in ( \tfrac{\ell}{L+\ell}, \tfrac L{L+\ell})$ for all~$t \in [0,T]$.
	\end{lemma}
	\begin{proof}  
	If $x(t)$ is identically constant then the proof follows from Lemma~\ref{l.constant}.
	Hence, we assume that~$x(t)$ is not constant. 
	This implies that we can restrict attention to the four regular cases depicted
	in  Figure~\ref{f.signs}.  
	
	Suppose that~$x(0) > \tfrac{L}{L+\ell}$.  Then 
	  considering the regular  cases 	depicted
	in  Figure~\ref{f.signs}, we see that~$x(T)<\tfrac{L}{L+\ell}$ and this contradicts the
	periodicity condition~$x(0)=x(T)$.

		Suppose that~$x(0) = \tfrac{L}{L+\ell}$.  Then 
	  considering the regular  cases 	depicted
	in  Figure~\ref{f.signs}, we see that again~$x(T)<  \tfrac{L}{L+\ell}$,
	as~$x(t)$ can increase towards~$  \tfrac{L}{L+\ell}$ only in the fourth case 
	 depicted
	in  Figure~\ref{f.signs}, yet it can never reach~$  \tfrac{L}{L+\ell}$, as this is an
	equilibrium   (and thus an invariant set) of this dynamics. 
	
	Summarizing, we showed that~$x(0)>\tfrac{\ell}{L+\ell}$. Using a  similar argument  shows  
	that~$x(0)\in (\tfrac{\ell}{L+\ell}  , \tfrac{L}{L+\ell} )$, and this implies that~$x(t) \in (\tfrac{\ell}{L+\ell}  , \tfrac{L}{L+\ell} )$ for all~$t\in[0,T]$. 
	   	Similar  arguments can be used to prove the corresponding statement for~$p $.
\end{proof}

	The following lemma excludes certain transitions between arcs.		

	\begin{lemma}\label{l.t1}
		Let $\mathscr X$ be an extremal trajectory. If there exists $\tau \in[0,T]$ such that  $\varphi_0(\tau) \varphi_1(\tau) < 0$, then  $\varphi_0(t) \varphi_1(t) < 0$ for all $t\in [\tau,T]$.
	\end{lemma}

\begin{figure}
		\centering
		\includegraphics[width=0.85\textwidth]{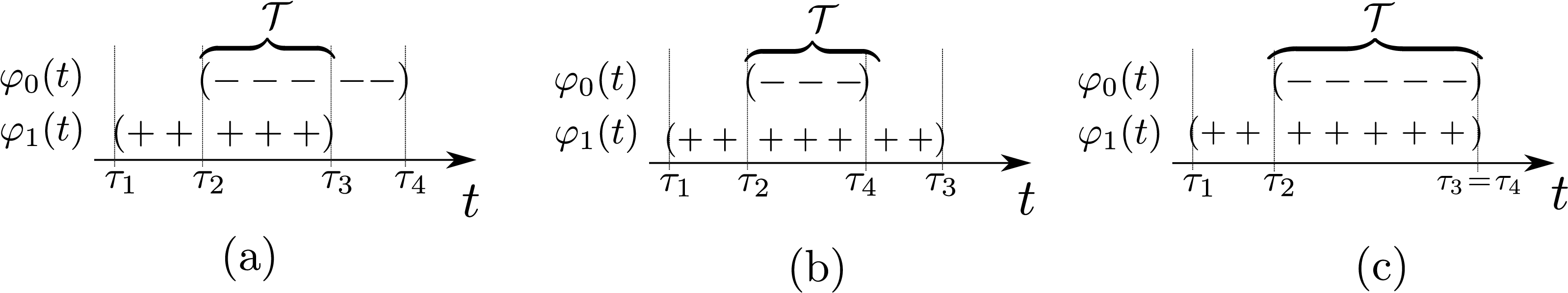}
		\caption{An illustration of the three cases studied in the proof of Lemma \ref{l.t1}: (a) $\tau_3 <\tau_4$, (b) $\tau_3>\tau_4$, (c) $\tau_3=\tau_4$.}\label{f.intervals}
	\end{figure}

	\begin{proof}
		W.l.o.g, assume that~$\varphi_0(\tau)<0$ and~$ \varphi_1(\tau)>0$.
		Hence,   $\tau \in E_{-}^0 \cap  E_{+}^1$.
		Since both sets are open,  there exists a connected component~$\cT \subset E_{-}^0 \cap E_{+}^1$ 
		such that~$\tau \in \cT$.
		Let $\cT_0, \cT_1$ be the connected components containing $\tau$ with respect to $E_-^0, E_+^1$, respectively. Then~$\cT = \cT_0 \cap \cT_1$. 
		Let~$\tau_i$, $i=1,\dots,4$, be such that~$\cT_1=(\tau_1,\tau_3)$, $ \cT_0=(\tau_2,\tau_4)$.  
		W.l.o.g, assume that~$\tau_1 \le \tau_2$.

		Assume first that $\tau_2>0$.  Then
	there are three possibilities:~$\tau_3<\tau_4$,
	$\tau_3>\tau_4$, and~$\tau_3=\tau_4$;  see
	Figure~\ref{f.intervals}. 
By  definition,  $\varphi_0(\tau_2)=0$, $\varphi_1(\tau_2)>0$.   Lemma \ref{l.bang} implies that~$u_0(t)=\ell, u_1(t)=L$ for all $t \in \cT$. By \eqref{dotphi0} and~\eqref{dotphi1}, both $\varphi_0$ and~$ \varphi_1$ are differentiable on $\cT$,
 and the right-derivative~$D_{\tau_2}^+ \varphi_0 $ exists.   Since $\varphi_0(\tau_2)=0$ and $\varphi_0(t)<0$ on $t\in \cT$, we have 
		\begin{align} \label{phi_ineq}
		0 &\geq  D_{\tau_2}^+ \varphi_0  \nonumber 
		\\&= u_1(\tau_{2}^+) (p(\tau_2)- (1-x(\tau_2)))\nonumber  \\&= L( p(\tau_2)+x(\tau_2)-1 ) .  
		\end{align}
		Recall that~$\dot x(t) = \ell - (\ell+L)  x(t)$, $\dot p =  (\ell+L) p(t) - L$ for $t\in \mathcal T$ (see Figure \ref{f.signs}), so
		\begin{align} \label{x_ineq}
		x(t) &= \left ( x(\tau_2) - \frac{\ell}{\ell+L} \right ) e^{-(\ell+L)(t-\tau_2)} + \frac {\ell}{\ell+L}\nonumber \\
		&< \left ( x(\tau_2) - \frac{\ell}{\ell+L} \right ) e^{(\ell+L)(t-\tau_2)} + \frac {\ell}{\ell+L},\\
		p(t) &= \left ( p(\tau_2) - \frac{L}{\ell+L} \right ) e^{(\ell+L)(t-\tau_2)} + \frac {L}{\ell+L},\nonumber
		\end{align}
		where the inequality \eqref{x_ineq} follows from the fact
		that~$x(t)> \frac{\ell}{\ell+L}$ (see Lemma \ref{l.bounds}).
Summing up these equations gives
\[
	x(t)+ p(t) <  ( x(\tau_2)+p(\tau_2) - 1 ) e^{(\ell+L)(t-\tau_2)} + 1.
	\]
		Thus, for any~$t \in \cT$, 
		\begin{align*}
		\dot\varphi_0(t)& = L ( p(t) + x(t) - 1 )\\
			&< L ( x(\tau_2)+p(\tau_2) -1  ) e^{(\ell+L)(t-\tau_2)}  \\& \le 0,
		\end{align*}
 		where the last  inequality follows from~\eqref{phi_ineq}. 
Hence, $\dot\varphi_0(t)<0$ on $\cT$.   Since $\sgn \dot\varphi_1(t) = -\sgn \dot \varphi_0(t)$,   $\dot\varphi_1(t)>0$ for all $t\in \cT$.

We now show that~$\tau_3=\tau_4$. 		
 Assume that $\tau_3 < \tau_4$. Then $\varphi_0(\tau_3)<0$ and $\varphi_1(\tau_3)=0$ as shown in Figure \ref{f.intervals}(a).   Integrating $\dot\varphi_1$ over $\cT$, and since $\varphi_1(\tau_2)\ge 0$, we get that 
		$\varphi_1(\tau_3)>0$, which is a contradiction. Similarly, assuming that~$\tau_3>\tau_4$ gives~$\varphi_0(\tau_4)=0$ (see Figure \ref{f.intervals}(b)). Since $\varphi_1(\tau_2)=0$ and $\dot\varphi_0(t)<0$ on $\cT$ then $\varphi_0(\tau_4)<0$, which is a contradiction.
Thus, $\tau_3=\tau_4$ (see Figure \ref{f.intervals}(c)).
		
		Let $\tau_e:=\tau_3=\tau_4$. Then the preceding argument shows   that~$\varphi_0(\tau_e)<0$ and~$\varphi_1(\tau_e)>0$.
		If~$\tau_e<T$ then the definition of~$\tau_3,\tau_4$ implies that~$\phi_0(\tau_e)\phi_1(\tau_e)=0$. 
		We conclude that~$\tau_e=T$. i.e. $\sup \cT = T$.  
		
		Assume now that $\inf \cT=0$.  Since $\mathscr X$ is periodic, we can study it on the interval $[0,2T]$. Let $\tilde E_-^0:=  E_-^0 \cup ( \{T\}+E_-^0) $, $\tilde E_+^1:=  E_+^1 \cup ( \{T\}+E_+^1)$.  Hence,  define $\tilde {\cT}$ as the maximal open neighborhood containing $\tau=T$ in $\tilde E_-^0 \cap \tilde E_+^1$. The sets  $ \tilde {\cT}_1, \tilde{\cT}_2$  are defined similarly.  Replicating the previous arguments to the sets $ \tilde{\cT}, \tilde {\cT}_1, \tilde{\cT}_2$ we see that $\sup \tilde{\cT} = 2T$. Hence, by periodicity, $\sup \cT = T$. 
	\end{proof}

	We can strengthen     Lemma~\ref{l.t1}  to exclude mixed-sign arcs.
	\begin{lemma}\label{l.no_mixed_sign}
		Let $\mathscr X$ be an extremal trajectory. Then,  $\varphi_0(t) \varphi_1(t) \ge 0   $ 
		for all~$t\in [0,T]$. 
	\end{lemma}

	\begin{proof}
		Assume that there exists $\tau \in [0,T]$ such that $\varphi_0(\tau)\varphi_1(\tau)<0$. 
			  Lemma \ref{l.t1}~implies that $\varphi_0(t)\varphi_1(t)<0$ for all $t \in [\tau,T]$.  By periodicity,~$\varphi_0(0)\varphi_1(0)<0$. Applying Lemma~\ref{l.t1} gives~$\varphi_0(t) \varphi_1(t) < 0   $ for all $t\in [0,T]$. 
				This implies that both $\varphi_0(t), \varphi_1(t)$ have constant and opposite signs.
W.l.o.g, assume that~$\varphi_0(t)<0$ and~$ \varphi_1(t)>0$ for all~$t \in [0,T]$. By Lemma~\ref{l.bang}, 
 $u_0(t)\equiv \ell$ and~$ u_1(t)\equiv  L$.   Hence,  
$\dot x(t)= \ell - (\ell + L ) x(t) $ for all~$t\in [0,T]$, so
		\[ x(T)-x(0) = \left (\frac{\ell}{\ell+L} - x(0) \right) (1- e^{-(\ell+L)T}). 
		\]
		Since  $x(0)>\frac{\ell}{\ell+L}$ (see Lemma \ref{l.bounds}),   $x(T)<x(0)$, and this is a contradiction.
	\end{proof}
	
	For an extremal trajectory $\mathscr X$, 
	recall that $E_r:=\{t  \in [0,T] \st  \varphi_0(t)\varphi_1(t) \ne 0\}$.   Lemma~\ref{l.no_mixed_sign} implies that
	\be\label{eq:prre} E_r= E_{++} \cup E_{--},
	\ee
	where
	\[
	E_{++}:=E_{+}^0 \cap E_+^1 \mbox{ and } 
	E_{--}:=E_{-}^0 \cap E_-^1.
	\]

In other words, the only possible   bang arcs are the first two cases in	Figure~\ref{f.signs}. Note that~$1/2$ is an equilibrium point of both these arcs. Also, on a singular arc~$x(t)$ is constant. This proves the following.

	\begin{lemma}\label{l.ge_half}
		Let $\mathscr X$ be an extremal trajectory. If $x(\tau)\ne 1/2$ for some~$\tau \in [0,T]$, then~$x(t) \ne 1/2$ for all~$t \in [0,T]$
	\end{lemma}

	The next lemma shows that an extremal trajectory must consist of a single singular arc.
	
		\begin{lemma}
		Let $\mathscr X$ be an extremal trajectory. Then,  $\varphi_0(t) \varphi_1(t)= 0   $ for all $t\in [0,T]$. Furthermore, $x (t)\equiv x(0)$ for all~$t\in [0,T]$.
	\end{lemma}
	\begin{proof}  
	We consider two cases.
	
	\noindent {\sl Case 1.} Suppose  that there exists a~$\tau\in[0,T]$ such that~$x(\tau)\ne 1/2$. 
We may assume w.l.o.g. that~$x(\tau)>1/2$. 
By Lemma~\ref{l.ge_half},
	$x(t)>1/2$ for all~$t\in[0,T]$.
	Seeking a contradiction, assume that~$\mu(E_r)>0$. Then~\eqref{eq:prre} implies that~$\dot x(t) <0$ for all~$t \in E_r$ (see Figure~\ref{f.signs}).
	By Lemma~\ref{l.s1},   $\dot x(t)=0$ for almost all $t\in E_s$.
	We conclude that~$x(T)<x(0)$, and this is a   contradiction. 
	Thus,~$\mu(E_r) = 0$.
	
	\noindent {\sl Case 2.} Suppose  that~$x(t)\equiv 1/2$.
	Then Lemma~\ref{l.constant} implies that~$\mu(E_r) = 0$.  
		\end{proof}

	We can now prove Theorem~\ref{th.extremal}.
 		We already know that~$\mathscr X$ consists of a  single singular arc. Lemma~\ref{l.s1} implies that there exists a~$c\in(0,1)$
		such that~$x(t)\equiv c$ for all~$t\in[0,T]$. 
		Integrating~\eqref{inputs_coupled} over~$[0,T]$ yields~$\left(\tfrac 1c -1\right) \bu_0 = \bu_1$, and this completes the proof.

\subsection*{A.3 \quad An alternative proof of Theorem \ref{theorem2in}}
  The alternative proof is inspired by the completing the square  idea in  \cite{katriel20} which was used to prove the result for a system with a controlled inflow and a constant outflow (proposed earlier by the authors in the preprint \cite{biorxiv}). We show that a similar and simpler approach can be developed to tackle 
  the more general case when both the inflow and outflow can be controlled independently of each other.  The proof is  based on two Lemmas. 
	Let $T>0$.  Recall that we use the notation  $\mean{y}:=\frac{1}{T}\int_0^T y(s)\,ds$. %
	\begin{lemma} \label{lem.s}
 Let~$x_1(t)$ be a solution of~\eqref{rfm1d} satisfying~$x_1(T)=x_1(0)$. Then
	$
	\mean{u_0 x_1^{k}} = \mean {(u_0+u_1) x_1^{k+1}} $
	 for any  integer $k\geq0$.
\end{lemma}
\begin{proof} 
Since $x_1^{k+1}(T) - x_1^{k+1}(0) = 0$,  
	the integral of
	\[
	\frac{1}{k+1}\frac{dx_1^{k+1}}{dt} = x_1^k \frac{dx_1}{dt}=x_1^k (u_0(1- x_1) -  u_1  x_1 )
	\]
	is zero, and the result follows.
\end{proof}

For~$k=0,1$, Lemma \ref{lem.s} gives
\begin{align}\label{eq:use6}
  {\bu_0}& = \overline {(u_0+u_1) x_1},\nonumber \\
\overline{u_0 x_1} &= \overline{(u_0+u_1) x_1^2},
\end{align}
 respectively.

\begin{lemma} \label{lem:pofgt}
 Let~$x_1(t)$ be a solution of~\eqref{rfm1d} satisfying~$x_1(T)=x_1(0)$.
Then
\be\label{eq:rzty}
\mean{u_1x_1}\leq    z,
\ee
where~$z:=   { \bu_0 }   {\bu_1} /( \bu_1 +\bu_0) $.
\end{lemma}
\begin{proof} 
Writing~$\mean{u_1x_1}=\mean{((u_0+u_1)-u_0)x_1} $,
and using~\eqref{eq:use6} gives~$\mean{u_1x_1}=\bu_0- \mean{u_0 x_1}
=\bu_0- \mean{(u_0+u_1)x_1^2} $.
To apply  a completion of squares argument, 
write  this as~$\mean{u_1x_1} = 
z  + \frac{ {\bu_0}    {\bu_0}}{ {\bu_0}+ \bu_1} 
- \mean{(u_0+u_1)x_1^2} $. Now~\eqref{eq:use6} gives
\begin{align*}
\mean{u_1x_1} 
  &   = z-
	\left ( \frac{ {\bu_0} \bu_0}{ \bu_0 + \bu_1} - \frac{2\bu_0} {\bu_0 +\bu_1} 
		\mean{ (u_0+u_1)x_1}  + \mean{(u_0+u_1)x_1^2}  \right ) \\
& = 
z-  \mean{(u_0+u_1)\left(x_1 - \frac{ {\bu_0}}{ {\bu_0}+ {\bu_1}}\right)^2},
\end{align*}
and this completes the proof.
\end{proof}

Eq.~\eqref{eq:rzty} implies that~$x_1(t)\equiv \frac{ {\bu_0}}{ {\bu_1}+ {\bu_0}}$ is an optimal
trajectory, thus providing an alternative proof to Theorem~\ref{th.extremal}.

		\subsection*{A.4 \quad Proof of Proposition \ref{pr.tv}}
The Hamiltonian is	$\mathcal H = p_1 ( u_0 (1-x_1)-\lambda_1 x_1) +p_2 u_0 +    \beta(t) \lambda_1 x_1 $, where we assume 
  w.l.o.g. that $p_0=T$. Hence, \eqref{pdot} gives 
\begin{align*}
 \dot p_1(t) &= (u_0(t) + \lambda_1 ) p_1(t) - \lambda_1 \beta(t),\\
\dot p_2(t)& \equiv 0.
\end{align*}
The switching function is
\begin{align}\label{eq:sewocp}
\varphi_0&= p_1 (1-x_1 )+p_2(0) , 
\end{align}
and thus
\begin{align}\label{eq:sewocpder}
\dot \varphi_0&= \lambda_1 (p_1-\beta +\beta x_1),\nonumber\\
\ddot \varphi_0&= \lambda_1 ((u_0 + \lambda_1 ) p_1 - \lambda_1 \beta -\dot \beta+
\dot \beta x_1 +\beta (u_0 (1-x_1)-\lambda_1 x_1) )  , 
\end{align}
  By Lemma \ref{l.bang}, an optimal control 
	is bang-bang on~$E_r$. Hence, we study the control on the set~$E_0:=\{t  \st  \varphi_0(t)= 0\}$ which we assume to have nonzero measure. As in  the proof of Lemma \ref{l.s1},  we can find a set~$F \subseteq E_0$ such that~$F=E_0$ a.e. and 
	$\varphi(t) = \dot \varphi(t)=\ddot \varphi(t) $  for all~$t\in F$.
  This gives
\begin{align}\label{eq:ubyd}
p_1(t) &= \beta(t) (1-x_1(t)),\nonumber \\
(1-x(t))^2 &=  -p_2(0) / \beta(t),\nonumber\\
u_0(t)&= \frac {\lambda_1 } {1-x(t)} - \lambda_1 +  \tfrac{\dot\beta(t)}{2\beta(t)},
\end{align}
and this proves~\eqref{e.tv_singular}. 
 Note that this implies that $p_2(0)<0$.  
 Furthermore, if $E_0=[0,1]$, then 
the equation $\int_0^1 u_0(t)   \diff t =  \bu_0$ yields~\eqref{constant}.

			\section*{Acknowledgments}
			This research has been  supported in part 
			by  research grants from    
			the Israel Science Foundation,  the 
			US-Israel Binational Science Foundation, and NSF 1817936. The authors thank Mahdiar Sadeghi for helpful discussions. 

\end{document}